\documentclass[11pt]{article}
\usepackage{amssymb}
\usepackage{graphicx}
\usepackage{amsmath}
\usepackage{amssymb}
\usepackage{graphicx}
\usepackage{amsmath}
\usepackage{amsthm}
\usepackage{enumerate}
\usepackage{paralist}
  
\usepackage[margin=1in]{geometry}
\def\inbar{\,\vrule height1.5ex width.4pt depth0pt}
\def\inbar{\,\vrule height1.5ex width.4pt depth0pt}
\def\IB{\relax{\rm I\kern-.18em B}}
\def\IC{\relax\hbox{$\inbar\kern-.3em{\rm C}$}}
\def\ID{\relax{\rm I\kern-.18em D}}
\def\IE{\relax{\rm I\kern-.18em E}}
\def\IF{\relax{\rm I\kern-.18em F}}
\def\IG{\relax\hbox{$\inbar\kern-.3em{\rm G}$}}
\def\IH{\relax{\rm I\kern-.18em H}}
\def\II{\relax{\rm I\kern-.18em I}}
\def\IK{\relax{\rm I\kern-.18em K}}
\def\IL{\relax{\rm I\kern-.18em L}}
\def\IM{\relax{\rm I\kern-.18em M}}
\def\IN{\relax{\rm I\kern-.18em N}}
\def\IO{\relax\hbox{$\inbar\kern-.3em{\rm O}$}}
\def\IP{\relax{\rm I\kern-.18em P}}
\def\IQ{\relax\hbox{$\inbar\kern-.3em{\rm Q}$}}
\def\IR{\relax{\rm I\kern-.18em R}}

\def\IZ{\relax\ifmmode\mathchoice
{\hbox{\cmss Z\kern-.4em Z}}{\hbox{\cmss Z\kern-.4em Z}}
{\lower.9pt\hbox{\cmsss Z\kern-.4em Z}}
{\lower1.2pt\hbox{\cmsss
Z\kern-.4em Z}}\else{\cmss Z\kern-.4em Z}\fi}
\def\IGa{\relax\hbox{${\rm I}\kern-.18em\Gamma$}}
\def\IPi{\relax\hbox{${\rm I}\kern-.18em\Pi$}}
\def\ITh{\relax\hbox{$\inbar\kern-.3em\Theta$}}
\def\IOm{\relax\hbox{$\inbar\kern-3.00pt\Omega$}}

\begin{document}

\baselineskip 20pt
\pagenumbering{arabic}
\pagestyle{plain}

\newtheorem{defi}{Definition}[section]
\newtheorem{theo}[defi]{Theorem}
\newtheorem{lemm}[defi]{Lemma}
\newtheorem{prop}[defi]{Proposition}
\newtheorem{note}[defi]{Note}
\newtheorem{nota}[defi]{Notation}
\newtheorem{exam}[defi]{Example}
\newtheorem{coro}[defi]{Corollary}
\newtheorem{rema}[defi]{Remark}
\newtheorem{cons}[defi]{Construction}
\newtheorem{ques}[defi]{Question}
\newtheorem{conj}[defi]{Conjecture}
\newtheorem{term}[defi]{Terminology}
\newtheorem{disc}[defi]{Discussion}
\newtheorem*{observation*}{Observation}

\newcommand{\abs}{${\bar A}^*}
\newcommand{\qbs}{${\bar Q}^*}
\newcommand{\be}{\begin{enumerate}}
\newcommand{\ee}{\end{enumerate}}
\newcommand{\fany}{\rm for\ \ any\ \ }

\newcommand{\rb}{\overline{R}}
\newcommand{\mt}{\overline{M}}
\newcommand{\nt}{\overline{N}}
\newcommand{\nb}{\widetilde{N}}
\newcommand{\mb}{\widetilde{M}}
\newcommand{\m}{\bf {m}}

\def\cm{Cohen-Macaulay}
\def\wrt{with respect to\ }
\def\pni{\par\noindent}
\def\wma{we may assume without loss of generality that\ }
\def\Wma{We may assume without loss of generality that\ }
\def\ets{it suffices to show that\ }
\def\bwoc{by way of contradiction}
\def\iff{if and only if\ }
\def\st{such that\ }
\def\fg{finitely generated}

\def\a{\goth a}
\def\p{\mathbf p}
\def\char{\rm char}
\def\isom{\thinspace \cong\thinspace}
\def\rtar{\rightarrow}
\def\rta{\rightarrow}
\def\l{\lambda}
\def\d{\Delta}

\title{\bf A NOTE ON ITOH (e)-VALUATION RINGS OF AN IDEAL}

\author{Youngsu Kim, Louis J.\ Ratliff, Jr.\, and
David E.\ Rush}

\maketitle

\begin{abstract}
Let  $I$  be a regular proper ideal in a Noetherian ring  $R$, let $e$ 
$\ge$ $2$  be an integer, let  $\mathbf T_e$ $=$ 
$R[u,tI,u^{\frac{1}{e}}]' \cap R[u^{\frac{1}{e}},t^{\frac{1}{e}}]$ 
(where  $t$  is an indeterminate and  $u$ $=$ $\frac{1}{t}$), and let  
$\mathbf r_e$ $=$ $u^{\frac{1}{e}} \mathbf T_e$.  Then the Itoh 
(e)-valuation rings of  $I$ are the rings $(\mathbf T_e/z)_{(p/z)}$, 
where  $p$  varies over the (height one) associated prime ideals of  
$\mathbf r_e$  and $z$  is the (unique) minimal prime ideal in  $\mathbf 
T_e$  that is contained in  $p$.  We show, among other things:

\noindent (1)
$\mathbf r_e$  is a radical ideal if and only if  $e$  is a common 
multiple of the Rees integers of  $I$.

\noindent
(2) For each integer  $k$ $\ge$ $2$, there is a one-to-one
correspondence
between the Itoh (k)-valuation rings $(V^*,N^*)$  of  $I$
and the Rees valuation
rings  $(W,Q)$  of  $uR[u,tI]$; namely, if  $F(u)$  is the quotient
field of  $W$, then $V^*$  is the integral closure of
$W$  in  $F(u^{\frac{1}{k}})$.

\noindent
(3)
For each integer  $k$ $\ge$ $2$, if  $(V^*,N^*)$  and  $(W,Q)$ are
corresponding
valuation rings, as in (2),
then  $V^*$  is a finite integral
extension domain of  $W$,
and  $W$  and  $V^*$  satisfy the Fundamental Equality with no 
splitting.
Also, if  $uW$ $=$ $Q^e$, and
if the greatest common divisor of  $e$  and  $k$  is  $d$,  and  $c$  is 
the
integer such that  $cd$ $=$ $k$,
then  $QV^*$ $=$ ${N^*}^c$ and $[(V^*/N^*):(W/Q)]$ $=$ $d$.
Further, if  $uW$ $=$ $Q^e$ and  $k$ $=$ $qe$  is a multiple of $e$,
then there exists a unit  $\theta_{e}$ $\in$ $V^*$  such that  $V^*$
$=$ $W[\theta_{e},u^{\frac{1}{k}}]$  is a finite free integral 
extension
domain of  $W$, $QV^*$ $=$ ${N^*}^q$,
$N^*$ $=$ $u^{\frac{1}{k}}V^*$, and $[V^*:W]$ $=$ $k$.

\noindent
(4)
If the Rees integers of  $I$  are all equal
to  $e$,
then  $V^*$ $=$ $W[\theta_e]$  is a simple free integral extension
domain of  $W$, $QV^*$ $=$ $N^*$ $=$ $u^{\frac{1}{e}}V^*$,
and  $[V^*:W]$ $=$ $e$ $=$ $[(V^*/N^*):(W/Q)]$.
\end{abstract}

\section{INTRODUCTION}
All rings in this paper are commutative and have an identity element
$1$ $\ne$ $0$, and our terminology is mainly as in Nagata \cite{N2}.
Thus a {\bf basis} of
an ideal is a generating set of the ideal, the term {\bf altitude}
refers to what is often called  dimension or Krull dimension, and
for a pair of local rings  $(R,M)$  and  $(S,N)$, $S$ {\bf dominates}
$R$
in case  $R$ $\subseteq$ $S$  and  $N \cap R$ $=$ $M$, and
we then write
$\mathbf{R \le S}$  or  $\mathbf{(S,N) \ge (R,M)}$.

In 1988, Shiroh Itoh proved the following interesting and
useful theorem in
\cite[p. 392, lines 3-11]{It} (the terminology is
defined in Section 2):

\begin{theo}
\label{IT}
Let  $I$  be a regular proper ideal in a Noetherian ring  $R$,
let  $\mathbf R$ $=$ $R[u,tI]$  be the
Rees ring of  $R$  with respect to  $I$,
let  $(W_1,Q_1),\dots,(W_n,Q_n)$  be the Rees valuation rings
of  $u \mathbf R$, for $j$ $=$ $1,\dots,n$, let $uW_j$ $=$ 
${Q_j}^{e_j}$
\textup{(}so  $e_1,\dots,e_n$  are the Rees integers of  $I$\textup{)},
and let  $e$ $\ge$ $2$  be an arbitrary common multiple
of  $e_1,\dots,e_n$.  Also, let  $\mathbf S$
$=$ $\mathbf R[u^{\frac{1}{e}}]$, let  $\mathbf T$
$=$ $\mathbf S' \cap R[u^{\frac{1}{e}},t^{\frac{1}{e}}]$,
and let  $\mathbf r$ $=$ $u^{\frac{1}{e}} \mathbf T$. Then:

\noindent
{\bf{(\ref{IT}.1)}}
$\mathbf r$ is a radical ideal, so the Rees integers of  $\mathbf r$
and of  $(u^{\frac{1}{e}} \mathbf S)_a$  are all equal to one.

\noindent
{\bf{(\ref{IT}.2)}}
There is a one-to-one correspondence
between the Rees valuation rings  $(V^*,N^*)$  of  $\mathbf r$  and
the Rees valuation rings  $(W,Q)$  of  $u \mathbf R$; namely,
if  $F(u)$  is the quotient field of  $W$, then  $V^*$  is the
integral closure of  $W$  in  $F(u^{\frac{1}{e}})$.

\noindent
{\bf{(\ref{IT}.3)}}
Let  $(V^*,N^*)$  and  $(W,Q)$  be
corresponding
Rees valuation
rings of  $\mathbf r$  and  $u \mathbf R$, respectively,
as in \textup{(}\ref{IT}.2\textup{)},
so  $W$ $=$ $W_j$  for some  $j$ $\in$ $\{1,\dots,n\}$.  Then  $QV^*$
= ${N^*}^{\frac{e}{e_j}}$,
so the ramification index of  $V^*$  relative to  $W$  is equal to  
$\frac{e}{e_j}$.

\end{theo}

Actually, the only part of this
theorem that S. Itoh specifically stated
in \cite{It}
was that  $\mathbf r$  is a radical ideal
when  $e$  is the least common multiple of
$e_1,\dots,e_n$.
His proof of this
essentially shows that (\ref{IT}.1) - (\ref{IT}.3) hold,
but his goals in \cite{It} were to prove
several nice applications of the radicality of the ideal  $\mathbf r$,
not to find additional properties of the Rees valuation rings
of this ideal.

However, it turns out that the Rees valuation rings of ideals like
$\mathbf r$
have some additional nice properties, and the goal of this present paper
is to derive some of these properties.
To facilitate discussing these valuation rings we make the following
definition.

\begin{defi}
\label{DD}
{\em
For an arbitrary
integer  $e$ $\ge$ $2$, the {\bf{Itoh (e)-valuation rings of  $I$}} are
the Rees valuation rings of  $\mathbf r_e$ $=$
$u^{\frac{1}{e}} \mathbf T_e$, where  $\mathbf T_e$ $=$
$ R[u,tI,u^{\frac{1}{e}}]' \cap R[u^{\frac{1}{e}}, t^{\frac{1}{e}}]$.
}
\end{defi}

\section{DEFINITIONS AND KNOWN RESULTS}

In this section we recall the needed definitions and
mention the needed known results concerning them.

\begin{defi}
\label{defi1}
{\em
Let  $I$  be an ideal in a ring  $R$.  Then:

\noindent
{\bf{(\ref{defi1}.1)}}
$R'$  denotes the {\bf{integral closure}} of $R$ in its total quotient 
ring.

\noindent
{\bf{(\ref{defi1}.2)}}
$I_a$  denotes the
{\bf{integral closure}}
of  $I$  in  $R$, so  $I_a$  is the ideal
$\{x \in R \mid x \text{  is a root of an} \allowbreak
\text{equation of the form } X^n + i_1X^{n-1} + \cdots +i_n=0 \}$,
where  $i_j \in I^j$  for  $j =1,\dots,n$.
The ideal $I$  is
{\bf{integrally closed}}
in case  $I$ $=$ $I_a$.

\noindent
{\bf{(\ref{defi1}.3)}}
The
{\bf{Rees ring of}}
$R$
{\bf{with resect to}}
$I$  is the
graded subring $\mathbf R(R,I)$ $=$ $R[u,tI]$  of
$R[u,t]$, where  $t$  is an indeterminate and  $u$
$=$ $\frac{1}{t}$.

\noindent
{\bf{(\ref{defi1}.4)}}
Assume that  $R$  is Noetherian and that  $I$  is
a regular proper ideal in  $R$
(that is, $I$ contains a regular element of  $R$  and  $I$ $\ne$ $R$), 
and
let
$b_1,\dots,b_g$  be
regular elements in  $R$  that generate  $I$.  For
$i$ $=$ $1,\dots,g$,  let  $C_i$ $=$ $R[I/b_i]$,
let $p_{i,1},\dots,p_{i,n_i}$  be the
(height one)
associated prime ideals of  $b_i{C_i}'$
(see (\ref{defi1}.1)), let  $z_{i,j}$  be the
(unique; see Remark \ref{RR}.1 below)
minimal prime ideal in
${C_i}'$  that is contained in  $p_{i,j}$
(possibly  $z_{i,h}$ $=$ $z_{i,j}$  for some
$h$ $\ne$ $j$  or  $z_{i,j}$ $=$ $(0)$), and let
$V_{i,j}$ $=$ $({C_i}'/z_{i,j})_{(p_{i,j}/z_{i,j})}$,
so  $V_{i,j}$ is a DVR.
Then the
set $\mathbf{RV}(I)$  of all  $V_{i,j}$  ($i$ $=$ $1,\dots,g$  and
$j$ $=$ $1,\dots,n_i$)  is the set of
{\bf{Rees valuation rings}}
of  $I$.  (The Rees valuation rings of  $I$
are well defined by  $I$; they do not depend
on the basis  $b_1,\dots,b_g$  of  $I$.)

\noindent
{\bf{(\ref{defi1}.5)}}
If  $(V,N)$  is a Rees valuation ring of  $I$
(see (\ref{defi1}.4)), then the
{\bf{Rees integer}}
of  $I$
{\bf{with respect to}}
$V$  is the positive integer  $e$  such that
$IV$ $=$ $N^e$.
If  $\mathbf{RV}(I)$ $=$ 
$\{(V_{1,1},N_{1,1}),\dots,(V_{g,n_g},N_{g,n_g})
\}$,
and if  $IV_{i,j}$ $=$
${N_{i,j}}^{e_{i,j}}$  for  $i$ = $1,\dots,g$
and  $j$ $=$ $1,\dots,n_i$,
then
$e_{1,1},\dots,e_{g,n_g}$  are called the
{\bf{Rees integers}}
of  $I$.

\noindent
{\bf{(\ref{defi1}.6)}}
If  $(W,Q)$ $\le$ $(V^*,N^*)$ are DVRs such that
$V^*$  is a localization of an integral extension domain of  $W$, then 
the
{\bf{ramification index of}}
$V^*$
{\bf{relative to}}
$W$  is the positive integer
$k$  such that  $QV^*$ $=$ ${N^*}^k$.
}
\end{defi}

\begin{rema}
\label{RR}
{\em
Let  $I$  be a regular proper ideal in a Noetherian ring  $R$.  Then:

\noindent
{\bf{(\ref{RR}.1)}}
Concerning (\ref{defi1}.4), it is shown in \cite[Definition p.\ 213 and
Propositions 2.7 and 2.13]{R2}
that if  $b$  is a regular nonunit in the integral closure  $A'$
of a Noetherian
ring $A$, then $bA'$  has a finite primary decomposition,
each associated prime ideal $p$ of $bA'$ has height one,
$p$  contains exactly one associated prime ideal  $z$  of  $(0)$,
and  $(A'/z)_{(p/z)}$  is a DVR.

\noindent
{\bf{(\ref{RR}.2)}}
It is shown in \cite[Proposition 10.2.3]{SH} that if  $V_1,\dots,V_n$  
are
the Rees valuation rings of  $I$, then  $(I^k)_a$ $=$
$\cap \{ I^k V_j \cap R \mid j = 1,\dots,n\}$  for all
$k$ $\in$ $\mathbb N_{> 0}$.  ($\mathbb N_{> 0}$  denotes the
set of positive integers.)

\noindent
{\bf{(\ref{RR}.3)}}
It readily follows from Definition \ref{defi1}.4 that the set
$\mathbf{RV}(I)$  of Rees valuation rings of  $I$  is the disjoint
union of the sets  $\mathbf{RV}((I+z)/z)$,
where  $z$  runs over the minimal prime ideals  $w$ in  $R$  such
that  $I+w$ $\ne$ $R$.  Also, if  $V$  is a
Rees valuation ring  of  $I$  and of  $(I+z)/z$, then the
Rees integer of  $I$  with respect to  $V$  is the
Rees integer of  $(I+z)/z$  with respect to  $V$.

\noindent
{\bf{(\ref{RR}.4)}}
It follows from Definition \ref{defi1}.4 (and
is proved in \cite[Example 10.3.2]{SH}) that if  $I$ $=$
$bR$  is a regular proper principal ideal in  $R$, then
the Rees valuation rings of  $I$  are the rings
$(R'/z)_{(p/z)}$, where  $p$  varies over the
(height one) associated prime ideals  of  $bR'$
and  $z$  is the (unique) minimal prime ideal in  $R'$  that
is contained in  $p$.  (It is readily checked that
if  $A$ $=$ $R' \cap R[\frac{1}{b}]$, then there is a
one-to-one correspondence between the associated
prime ideals  $p$  of  $bR'$  and the associated
prime ideals  $q$ of  $bA$; namely, $q$ $=$ $p \cap A$,
and then  $(R'/z)_{(p/z)}$
$=$ $(A/z')_{(q/z')}$, where  $z'$ $=$ $z \cap A$.)

\noindent
{\bf{(\ref{RR}.5)}}
It follows from Definitions \ref{defi1}.4 and
\ref{defi1}.5 that if
$A$  is a ring such that  $R$ $\subseteq$ $A$ $\subseteq$
$R'$, then  $I$  and  $IA$  have the same Rees valuation
rings and the same Rees integers.  And it also follows
that, for all positive integers  $k$,
$IA$  and  $I^kA$  have the same Rees valuation rings.

\noindent
{\bf{(\ref{RR}.6)}}
It is well known (and is
readily proved, much as in the proof of
\cite[(2.5)]{CHRR3}), that if  $R$  is a Noetherian integral
domain and  $A$  is an integral extension domain of  $R$, then
the Rees valuation rings of  $IA$  are
the extensions of the Rees valuation rings of  $I$
to the quotient field of  $A$.
}

\end{rema}

The following theorem is a special case of
\cite[Theorems 19 and 20, pp. 55 and 60-61]{ZS2};
the terminology ``Fundamental Inequality'' is due to
O. Endler in \cite[pp. 127-128]{End}.

\begin{theo}
\label{FI}
\textup{(}Fundamental Inequality\textup{)}: Let $(V,N)$  be a DVR with
quotient
field  $F$, let  $E$  be a finite algebraic extension field of
$F$, let  $[E:F]$ $=$ $e$, let  $(V_1,N_1),\dots,(V_n,N_n)$
be all of the valuation rings of  $E$  that are extensions
of  $V$  to  $E$  \textup{(}so the integral closure  $V'$  of  $V$  in
$E$  has exactly  $n$  maximal ideals
$M_1,\dots,M_n$  and  $V_j$ $=$ ${V'}_{M_j}$, $j$
$=$ $1,\dots,n$\textup{)}, and for
$j$ $=$ $1,\dots,n$,  let  $NV_j$ $=$ ${N_j}^{e_j}$  and
$[(V_j/N_j):(V/N)]$ $=$ $f_j$.  Then
$$\textup{(FI)}~~~~~~~~~~\overset{n}{\sum_{j=1}} \; e_jf_j \le e,$$
and the equality holds if the integral closure  $V'$  of
$V$  in  $E$  is a finite  $V$-module.

\end{theo}

\begin{term}
\label{term1}
{\em
If  $(V,N)$  and  $(V_1,N_1),\dots,(V_n,N_n)$  are as in
Theorem \ref{FI}, if  $n$ $=$ $1$, and if the equality in (FI) holds,
then
it will be said that
$(V,N)$  and  $(V_1,N_1)$
{\bf{satisfy the Fundamental Equality with no splitting}}.
}

\end{term}

\begin{nota}
\label{nota1}
{\em
If  $R$  is an integral domain, then  $R_{(0)}$  denotes the
quotient field of  $R$.  Therefore, if  $S$  is a finite algebraic
extension domain of  $R$, then  $[S_{(0)} : R_{(0)}]$  denotes
the dimension of the quotient field  $S_{(0)}$  of  $S$  over the
quotient field  $R_{(0)}$  of  $R$.
If  $S$  is a finite free integral extension domain of  $R$,
and if it is clear that the rank of  $S$  is equal to
$[S_{(0)} : R_{(0)}]$, then we often write  $[S:R]$  in place of
$[S_{(0)} : R_{(0)}]$.
}

\end{nota}

The next three propositions are known, but we do not
know specific references for them, so we sketch their proofs.

\begin{prop}
\label{prop1}
Let  $(W,Q)$  be a  DVR, let  $Q$ $=$ $\pi W$,
let  $f$ $\ge$ $2$  be an
integer, let  $D$ $=$ $W[\pi^{\frac{1}{f}}]$,
and let  $P'$ $=$ $\pi^{\frac{1}{f}}D$.
Then  $D$  is a DVR that is a simple free integral extension domain of
$W$, $P'$  is the maximal ideal in  $D$,
$QD$ $=$ ${P'}^f$,
$[D_{(0)}:W_{(0)}]$ $=$ $f$, and  $D/P'$ $\cong$ $W/Q$.
Therefore  $W$  and  $D$  satisfy the Fundamental
Equality with no splitting \textup{(}see Terminology 
\ref{term1}\textup{)}.

\end{prop}

\begin{proof}
${P'}$ $=$ $(Q,\pi^{\frac{1}{f}})D$  is a maximal ideal
in  $D$, since  $D/{P'}$ $\cong$ $W/Q$  is a field, and  ${P'}$
$=$ $\pi^{\frac{1}{f}}D$  is a principal ideal,
since  $Q$  is generated
by  $\pi$.  Also, $D$  is integral over  $W$  and  $Q$  is
the unique maximal ideal in  $W$, so every maximal ideal
in  $D$  contains $QD$, so it follows that
${P'}$  is the only maximal ideal in  $D$,
hence  $D$  is a DVR.  It therefore follows that
$QD$ $=$ ${P'}^f$.  Therefore, since  $[D_{(0)}:W_{(0)}]$ $\le$ $f$,
it follows from the Fundamental Inequality
{FI} that  $[D_{(0)}:W_{(0)}]$ $=$ $f$, hence  $D$  is a simple free
integral extension domain of  $W$.
Thus, by (\ref{term1}), the last statement is clear.
\end{proof}

\begin{prop}
\label{prop2}
Let  $M$   be a maximal ideal in a Noetherian ring  $R$,
and let  $m(X)$  be a monic polynomial in  $R[X]$.
If the image  $\overline m (X)$  of  $m(X)$  in
$(R/M)[X]$  is irreducible, then  $m(X)$  is
irreducible in  $R[X]$, $MR[x]$  is a
maximal ideal in  $R[x]$ $=$
$R[X]/(m(X)R[X])$, and  $R[x]$  is a simple
free  integral extension ring of  $R$   of rank equal to  $\deg(m(X))$.

\end{prop}

\begin{proof}
By considering the maps
$R[X]$ $\rightarrow$
$(R/M)[X]$
$\rightarrow$ $(R/M)[\chi ]$,
where  $\chi$  is a root of the irreducible
(by hypothesis) polynomial  $\overline m(X)$  in
$(R/M)[X]$,
it follows that:
$m(X)$  is irreducible in $R[X]$;
$(M,m(X))R[X]$  is a maximal ideal;
$MR[x]$ $=$ $((M,m(X))R[X])/(m(X)R[X])$
is a maximal ideal in  $R[x]$; and,
$R[x]$  is a simple free integral extension
ring of  $R$  of rank equal to  $\deg(m(X))$.
\end{proof}

\begin{prop}
\label{prop3}
Let  $I$  be a regular proper
ideal in a Noetherian ring  $R$, let  $b_1,\dots,b_g$  be
regular elements in  $I$  that are a basis of  $I$,
and let  $\mathbf R$ $=$
$R[u,tI]$  be the Rees ring of  $R$  with respect to  $I$.
Then there exists a one-to-one correspondence between the
Rees valuation rings  $(V,N)$  of  $I$  and the Rees
valuation rings  $(W,Q)$  of  $uR[u,tI]$.  Namely,
if  $V$ $=$ $(R[I/b_i]'/z)_{(p/z)}$
\textup{(}where  $i$ $\in$
$\{1,\dots,g\}$, $p$  is an associated
prime ideal of  $b_iR[I/b_i]'$, and  $z$  is the
minimal prime ideal in  $R[I/b_i]'$  that is
contained in  $p$\textup{)}, then
$W$ $=$ $V[t \overline{b_i}]_{NV[t \overline{b_i}]}$,
where  $t \overline{b_i}$ $=$ $t(b_i + z)$,
and  $Q$ $=$ $NW$.  \textup{(}Note that
$\overline{tb_i}$ $=$ $tb_i + z'$
\textup{(}where  $z'$ $=$ $zR[u,t] \cap R[u,tI]$\textup{)}
corresponds to $t \cdot \overline {b_i}$  in
the isomorphism between  $R[u,tI]/z'$  and  $(R/z)[u,t((I+z)/z)]$;
see  \textup{\cite{ReF}}.\textup{)}
\end{prop}

\begin{proof}
Since  $t$  is transcendental over  $R$  and  $u$ $=$
$\frac{1}{t}$, there is a one-to-one correspondence
between the minimal prime ideals  $z$  in  $R$  and
the minimal prime ideals  $z'$  in  $\mathbf R$; namely,
$z'$ $=$ $zR[u,t] \cap \mathbf R$,
and then  $z$ $=$ $z' \cap R$.  Thus it follows from
Remark \ref{RR}.3 that it suffices to prove this
proposition for the case when  $R$  is a Noetherian
integral domain.

Therefore assume that  $R$  is a Noetherian domain,
fix  $b$ $\in$ $\{b_1,\dots,b_g\}$, let  $C$ $=$
$R[I/b]$, and let  $\mathbf C'$ $=$ $\mathbf R[\frac{1}{tb}]'$,
so  $\mathbf R'[\frac{1}{tb}]$ $=$
$\mathbf C'$ $=$ $C'[tb,\frac{1}{tb}]$.
Also, $u \mathbf C'$ $=$ $b \mathbf C'$, so, since
$tb$  is transcendental over  $C'$, there exists
a one-to-one correspondence between the (height one)
associated prime ideals  $p$  of  $bC'$  and the
(height one) associated prime ideals  $p'$  of  $u \mathbf C'$;
namely, $p'$ $=$ $p \mathbf C'$  (and  $p$ $=$
$p' \cap C'$).  Therefore there exists a one-to-one
correspondence between the DVRs  $V$ $=$ ${C'}_p$  and
the DVRs  $W$ $=$ ${\mathbf C'}_{p'}$, and for
corresponding  $V$ and  $W$, if $V$ $=$ ${R[I/b]'}_p$  and
$N$ $=$ $pV$, then $W$ $=$ $V[tb]_{NV[tb]}$  and
$Q$ $=$ $NW$ $=$ $pW$.

It follows that, for each Rees valuation ring  $(V,N)$  of
$I$  of the form  $V$ $=$ ${R[I/b]'}_p$, the ring  $W$ $=$
$V[tb]_{pV[tb]}$  is a Rees valuation ring of  $u \mathbf R$,
and  $W \cap R_{(0)}$ $=$ $V$.

Finally, let  $W$  be a Rees valuation ring of  $u \mathbf R$,
say  $W$ $=$ ${\mathbf R'}_{p'}$, where  $p'$  is a
(height one) associated prime ideal of  $u \mathbf R'$
(see Remark \ref{RR}.4).  To complete the proof of the
one-to-one correspondence, it suffices to show that
there exists  $b$ $\in$ $\{b_1,\dots,b_g\}$  and a
Rees valuation ring  $(V,N)$  of  $I$  such that
$V$ $=$ ${R[I/b]'}_p$  and
$W$ $=$ $V[tb]_{NV[tb]}$.

For this, $tb$ $\notin$ $p'$  for
some  $b$ $\in$ $\{b_1,\dots,b_g\}$, by \cite[Lemma 3.2]{Rt}
(the assumption in \cite{Rt} that  $R$  is analytically
unramified is not used in the proof of Lemma 3.2).
Therefore  $W$ $=$ ${\mathbf R[\frac{1}{tb}]'}_{p''}$, where
$p''$ $=$ $p' W \cap \mathbf R[\frac{1}{tb}]'$  (so
$p''$  is a (height one) associated prime ideal of
$u \mathbf R[\frac{1}{tb}]'$ $=$
$b \mathbf R[\frac{1}{tb}]'$, and
$\mathbf R[\frac{1}{tb}]'$ $=$ $C'[tb,\frac{1}{tb}]$, where
$C$ $=$ $R[I/b]$.  Let  $p$ $=$ $p'' \cap C'$.
Then  $p$  is a (height one) associated prime
ideal of  $bC'$  and  $p''$ $=$ $pC'[tb,\frac{1}{tb}]$,
so it follows from the second preceding paragraph that $W$
$=$ $V[tb]_{NV[tb]}$, where  $V$ $=$ ${C'}_p$  and
$N$ $=$ $pV$.
\end{proof}

\section{PROPERTIES OF ITOH (e)-VALUATION RINGS}
In this section we show that Itoh (e)-valuation rings have
several nice properties.  For this, we
need the following proposition, which is essentially
a corollary of Proposition \ref{prop2}.

\begin{prop}
\label{prop4}
Let  $I$  be a regular proper ideal in a Noetherian ring  $R$,
let  $b_1,\dots,b_g$  be regular elements in  $I$  that generate
$I$, and
let  $\mathbf R$ $=$ $R[u,tI]$  be the Rees ring of  $R$
with respect to  $I$.  Let
$(V,N)$  and  $(W,Q)$  be corresponding \textup{(}as in Proposition
\ref{prop3}\textup{)} Rees valuation rings
of  $I$  and  $u \mathbf R$, respectively, say  $V$
$=$ $(R[I/b]'/z)_{(p/z)}$, where  $b$ $\in$ $\{b_1,\dots,b_g\}$,
$p$  is a \textup{(}height one\textup{)}
associated prime ideal of  $bR[I/b]'$, and  $z$  is
the minimal prime ideal in  $R[I/b]'$  that is contained
in  $p$, so  $W$ $=$ $V[t \overline{b}]_{NV[t \overline{b}]}$,
where  $\overline b$ $=$ $b+z$.  Let  $e$ $\ge$ $2$ be an
integer, let  $v$ $\in$ $V - N$,
let  $m_e(X)$ $=$ $X^e - \frac{v}{t \overline{b}}$, and let
$\theta_e$  be a root of  $m_e(X)$ in a fixed
algebraic closure of the quotient field of  $W$.  Then
$U$ $=$ $W[\theta_e]$  is a DVR that is a
simple free  integral extension
domain of  $W$, $QU$ is its maximal ideal, and
$[U:W]$ $=$ $e$ $=$ $[U/(QU): W/Q]$.
Therefore  $W$  and  $U$  satisfy the Fundamental
Equality with no splitting \textup{(}see Terminology 
\ref{term1}\textup{)}.

\end{prop}

\begin{proof}
Since  $W =V[tb]_{NV[tb]}$
(with  $tb$  transcendental over
$V$)  and  $Q$ $=$ $NW$  is an extended ideal from  $V$,
$t \overline{b} + Q$ $\in$ $W/Q$ is transcendental over
$V/N$, so it follows that the image  $\overline{m_e}(X)$
$\in$ $(W/Q)[X]$  of  $m_e(X)$  is irreducible, hence
$[(W/Q)[y]:(W/Q)]$ $=$ $e$, where  $y$ $=$
$X+ \overline{m_e}(X)(W/Q)[X]$.  Therefore
it follows from Proposition \ref{prop2} that:
$m_e(X)$  is irreducible in the UFD  $W[X]$;
$M$ $=$ $QU$ is a (principal) maximal ideal,
where  $U$ $=$ $W[\theta_e]$;
and, $[U:W]$ $=$ $e$.
Further, $U$  is integral over  $W$  and  $Q$  is
the unique maximal ideal in  $W$, so every
maximal ideal in  $U$  contains  the maximal ideal $QU$, hence  $U$
is a DVR that is a simple free  integral extension
domain of  $W$  and  $M$  is its maximal ideal.
Therefore, by (\ref{term1}), the last statement is clear.
\end{proof}

Our first theorem, Theorem \ref{theo1y}, is an expanded
version of Itoh's Theorem (see Theorem \ref{IT};
note that (\ref{IT}.3) is proved in (\ref{theo1y}.5)).

\begin{theo}
\label{theo1y}
Let  $I$  be a regular proper ideal in a Noetherian ring  $R$,
let  $b_1,\dots,b_g$  be regular elements in  $I$  that are a
basis of  $I$,
let  $(V_1,N_1),\dots,(V_n,N_n)$  be the Rees valuation rings of  $I$,
let  $e_1,\dots,e_n$  be the Rees integers of  $I$,
let  $e$ be an arbitrary common multiple of
$e_1,\dots,e_n$, and let  $f_j$ $=$ $\frac{e}{e_j}$  \textup{(}$j$
$=$ $1,\dots,n$\textup{)}.  Also,
let  $\mathbf R$ $=$ $R[u,tI]$  be the
Rees ring of  $R$  with respect to  $I$, let  $\mathbf S$
$=$ $\mathbf R[u^{\frac{1}{e}}]$, let  $\mathbf T$
$=$ ${\mathbf {S}}' \cap R[u^{\frac{1}{e}},t^{\frac{1}{e}}]$,
and let  $\mathbf r$ $=$ $u^{\frac{1}{e}} \mathbf T$.  Then:

\noindent
{\bf{(\ref{theo1y}.1)}}
$\mathbf r$ is a radical ideal.

\noindent
{\bf{(\ref{theo1y}.2)}}
There is a one-to-one correspondence
between the Itoh \textup{(}e\textup{)}-valuation rings $(V^*,N^*)$ of
$I$  and
the Rees valuation rings  $(W,Q)$  of  $u \mathbf R$; namely,
given  $W$, if  $F(u)$  is the quotient field of  $W$,
then  $V^*$  is the integral closure of  $W$  in  $F(u^{\frac{1}{e}})$.

\noindent
{\bf{(\ref{theo1y}.3)}}
Let  $(W,Q)$ and  $(V^*,N^*)$  be corresponding
\textup{(}as in \textup{(}\ref{theo1y}.2\textup{))} valuation rings,
and let  $(V,N = \pi V)$  be the Rees
valuation ring of  $I$  that corresponds
\textup{(}as in Proposition \ref{prop3}\textup{)} to  $W$,
so  $V$ $\le$ $W$ $\le$ $V^*$  and
$NW$ $=$ $Q$.  Assume that  $V$ $=$ $V_j$, so  $IV$ $=$ 
${N}^{e_j}$
\textup{(}$=$ ${N_j}^{e_j}$\textup{)}.

Then there exists a unit  $\theta$ $\in$ $V^*$  such that  $(U,P)$
$=$ $(W[\theta],N^* \cap W[\theta])$  is an Itoh $(e_j)$-valuation
ring of  $I$  that is a simple
free integral extension domain of  $W$  and  $P$ $=$ $QU$
is the maximal ideal in  $U$, so
the ramification index of  $U$  relative to  $W$  is equal to one
\textup{(}see Definition \ref{defi1}.6\textup{)}.

Also,  $QU$ $=$ $\pi U$
$=$ $yU$, where  $y$ $=$ $u^{\frac{1}{e_j}}$,
and  $[U:W]$ $=$ $e_j$  $=$ $[(U/P):(W/Q)]$.
Therefore  $W$  and  $U$  satisfy the Fundamental Equality
with no splitting.

\noindent
{\bf{(\ref{theo1y}.4)}}
Let  $(V,N)$ $\le$ $(U,P)$ $\le$ $(V^*,N^*)$  \textup{(}with
$V$ $=$ $V_j$\textup{)}
be as in \textup{(}\ref{theo1y}.3\textup{)}.
Then there exists a nonunit  $x$  $\in$  $V^*$
such that  $V^*$ $=$ $U[x]$  is a simple free
integral extension domain of  $U$, $N^*$ $=$
$xV^*$ \textup{(}where  $x$ $=$ $y^{\frac{1}{f_j}}$ with $y$ $=$
$u^{\frac{1}{e_j}}$
\textup{(}as in \textup{(}\ref{theo1y}.3\textup{))}, so $x$ $=$
$u^{\frac {1}{e}}$\textup{)}, and  $PV^*$
$=$ ${N^*}^{f_j}$,
so the ramification index of  $V^*$  relative to  $U$  is equal to 
$f_j$.
Also, $[V^*:U]$ $=$ $f_j$, and
$V^*/N^*$ $\cong$ $U/P$,
so  $U$  and  $V^*$  satisfy the Fundamental Equality
with no splitting.

\noindent
{\bf{(\ref{theo1y}.5)}}
$V^*$ $=$ $W[\theta,x]$, $N^*$ $=$
$xV^*$ $=$ $u^{\frac{1}{e}}V^*$, where  $\theta$  is as
in \textup{(}\ref{theo1y}.3\textup{)} and
$x$  is as in \textup{(}\ref{theo1y}.4\textup{)}, and
$V^*$  is a finite free integral extension
domain of  $W$.  Also, the ramification index
of  $V^*$  relative to  $W$  is equal to $f_j$,
$[(V^*/N^*):(W/Q)]$ $=$ $e_j$, and  $[V^*: W]$ $=$ $e$, so
$W$  and  $V^*$  satisfy the Fundamental Equality
with no splitting.

\noindent
{\bf{(\ref{theo1y}.6)}}
Assume that  $e_1$ $=$ $\cdots$ $=$ $e_n$ $=$ $e$,
and let  $(W,Q)$  and  $U$ $=$ $W[\theta]$  be as
in \textup{(}\ref{theo1y}.3\textup{)}.
Then $V^*$ $=$ $U$  is a simple free integral extension domain
of  $W$, $P$ $=$ $QU$ $=$  $u^{\frac{1}{e}}U$ is the maximal ideal 
in
$U$,
and  $[U:W]$ $=$ $e_j$  $=$ $[(U/P):(W/Q)]$.
Therefore  $W$  and  $V^*$ \textup{(}$=$ $U$\textup{)} satisfy
the Fundamental Equality
with no splitting, and  $u^{\frac{1}{e}} \mathbf T$
is a radical ideal, where  $\mathbf T$ $=$ $R[u,tI,u^{\frac{1}{e}} ]' 
\cap R[u^{\frac{1}{e}},t^{\frac{1}{e}}]$.
\end{theo}

\begin{proof}
We first prove (\ref{theo1y}.3) - (\ref{theo1y}.5).  For this, fix
an Itoh (e)-valuation ring  $(V^*,N^*)$  of  $I$.  Then,
by Definitions \ref{DD} and \ref{defi1}.4, there
exists a (height one) associated prime
ideal  $p^*$  of  $u^{\frac{1}{e}} \mathbf T$  and a
minimal prime ideal  $z^*$ $\subset$ $p^*$  such that
$V^*$ $=$ $(\mathbf T/z^*)_{(p^*/z^*)}$.  Since  $\mathbf T$
$\subset$ $R[u^{\frac{1}{e}},t^{\frac{1}{e}}]$  and
$\mathbf T$  and $R[u^{\frac{1}{e}},t^{\frac{1}{e}}]$ have the same
total quotient ring, it follows
that  $z$ $=$ $z^*R[u^{\frac{1}{e}},t^{\frac{1}{e}}] \cap R$
is a minimal prime ideal in  $R$  and  $z^*$ $=$
$zR[u^{\frac{1}{e}},t^{\frac{1}{e}}] \cap \mathbf T$.
Also, $z'$ $=$ $zR[u,t] \cap \mathbf R$  is a minimal prime
ideal in  $\mathbf R$  and  $z^* \cap \mathbf R$ $=$ $z'$.

Let  $F$  be the quotient field of  $R/z$, so  $F(u)$
(resp., $F(u^{\frac{1}{e}})$)  is the quotient field of
$\mathbf R/z'$  (resp., $V^*$  and  $\mathbf T/z^*$).
Also, by Definition \ref{DD}, $V^*$  is a Rees valuation ring of
$u^{\frac{1}{e}} \mathbf T$, so  $V^*$  is also a Rees valuation
ring of $(u^{\frac{1}{e}} \mathbf T + z^*)/z^*$, by Remark \ref{RR}.3,
so  $V^*$  is also a Rees valuation ring of $(u \mathbf T + z^*)/z^*$,
by Remark \ref{RR}.5.

Therefore  $W$ $=$ $V^* \cap F(u)$  is a Rees valuation ring of the 
ideal
$(u \mathbf R + z)/z$, by Remark \ref{RR}.6 (and so
also of  $u \mathbf R$, by Remark \ref{RR}.3),
hence  $V^*$  is an extension of  $W$ to  $F(u^{\frac{1}{e}})$.

Let $V$  be the Rees valuation ring of  $I$  that
corresponds (as in Proposition \ref{prop3}) to  $W$,
and assume that  $V$ $=$
$V_j$ $=$ $(C'/z)_{(p/z)}$, where  $C'$
is the integral closure of  $C$ $=$ $R[I/b]$
in its total quotient ring,
$b$ $\in$ $\{b_1,\dots,b_g\}$,
$p$  is a (height one) associated
prime ideal of  $bC'$, and  $z$  is the (unique) minimal
prime ideal in  $C'$  that is contained in  $p$.
Let  $N$  and  $Q$  be the maximal ideals in  $V$  and  $W$,
respectively.  Then  $Q$ $=$ $NW$ by Proposition \ref{prop3}, so
\begin{displaymath}
(\ref{theo1y}.7) \; \; \;  t \overline{b} +Q \in W/Q
\mbox{ is transcendental over } V/N, \mbox{ where } \overline b = b+z.
\end{displaymath}

Let $N$ $=$ $\pi V$, so, by hypothesis,  $\overline{b}V$ $=$ $IV$ 
$=$
$N^{e_j}$$=$ $\pi^{e_j} V$.
Therefore  $uW$ $=$ $\overline{b}W$ $=$ $IW$ $=$ $Q^{e_j}$
$=$ $\pi^{e_j}W$, so
\begin{displaymath}
(\ref{theo1y}.8) \mbox{ there exist units } v \in V \mbox{ and } w \in W
\mbox{ such that } \overline{b} = v \pi^{e_j} \mbox{ and } u = w 
\pi^{e_j}.
\end{displaymath}

Since  $\overline{b}$ $=$ $u(t \overline{b})$
(see the last part of Proposition \ref{prop3}), it follows from
(\ref{theo1y}.8) that
$$(\ref{theo1y}.9)~~~~~~~~~~w = \frac{v}{t \overline{b}}.$$
Let $\theta$ $=$ $w^{\frac{1}{{e_j}}}$
(in a fixed algebraic closure  $(F(u))^*$  of  $F(u)$),
and let  $U$ $=$ $W[\theta]$.
Then it follows from Proposition \ref{prop4} that:
$U$  is a DVR that is a
simple free  integral extension domain of  $W$; $QU$
is its maximal ideal;
$[U:W]$ $=$ $e_j$ $=$ $[(U/(QU)):(W/Q)]$;
and, $W$  and  $U$  satisfy the Fundamental Equality with no splitting.

Continuing with the proof of (\ref{theo1y}.3), we next show that
$U$ $\le$ $V^*$, and we first show that
$\theta$ $\in$ $V^*$.  For this,
$w \pi^{e_j}$ $=$ $u$, by (\ref{theo1y}.8), $\theta$ $=$
$w^{\frac{1}{e_j}}$,
and  $u^{\frac{1}{e}}$ $\in$ $V^*$, so  $\theta \pi$ $=$
$(w \pi^{e_j})^{\frac{1}{e_j}}$ $=$ $u^{\frac{1}{e_j}}$ $=$ 
$(u^{\frac{
1}{e}})^{f_j}$
$\in$ $V^*$, hence  $\theta \pi$  $=$  $u^{\frac{1}{e_j}}$
$\in$ $V^*$.  Further, $\pi$ $\in$ $Q$ $\subset$ $V^*$, so  $\theta$  is
in the quotient field  $F(u^{\frac{1}{e}})$
of  $V^*$.  Moreover,  $\theta$ $=$ $w^{\frac{1}{e_j}}$  is
integral over  $W$  and is a unit, $W$ $\le$ $V^*$, and  $V^*$  is 
integrally closed in  $F(u^{\frac{1}{e}})$,
so  $\theta$ $\in$ $V^*$  and is a unit in  $V^*$.
Therefore  $U$ $=$ $W[\theta]$ $\le$ $V^*$,
and $u^{\frac{1}{e_j}}U$ $=$ $\theta \pi U$
$=$ $\pi U$ $=$ $QU$ $=$ $N^* \cap U$
is the maximal ideal in  $U$.
Thus, to complete the proof of (\ref{theo1y}.3),
it remains to show that  $U$  is
an Itoh ($e_j$)-valuation ring of  $I$.

For this, $U$ $=$ $W[\theta]$ $\supseteq$ $W[u^{\frac{1}{e_j}}]'$
$\supseteq$ $\mathbf R[u^{\frac{1}{e_j}}]'$ $\supseteq$
$\mathbf T_{e_j} + \mathbf R'$,
where  $\mathbf T_{e_j}$ $=$
$ R[u,tI,u^{\frac{1}{e_j}}]' \cap 
R[u^{\frac{1}{e_j}},t^{\frac{1}{e_j}}]$,
and  $P \cap \mathbf R'$ $=$
$(P \cap W) \cap \mathbf R'$ $=$ $Q \cap \mathbf R'$  is
a height one associated prime ideal of  $u \mathbf R'$,
so  $P \cap \mathbf T_{e_j}$ $=$ $q$, say,
is a height one associated prime
ideal of  $u^{\frac{1}{e_j}} \mathbf T_{e_j}$  (since
$u^{\frac{1}{e_j}}$ $\in$ $P$).

Therefore  $U$ $\ge$ $(\mathbf T_{e_j})_q$, $(\mathbf T_{e_j})_q$  is
an Itoh ($e_j$)-valuation ring of  $I$, and
$U$  and $(\mathbf T_{e_j})_q$  are DVRs with the same quotient
field, so  $U$ $=$ $(\mathbf T_{e_j})_q$  is an Itoh ($e_j$)-valuation
ring of  $I$.  Thus (\ref{theo1y}.3) holds.

For (\ref{theo1y}.4), let
$y$ $=$ $u^{\frac{1}{e_j}}$, let  $x$ $=$ $y^{\frac{1}{f_j}}$
(so  $x$ $=$
$u^{\frac{1}{e}}$, since  $e_j \cdot f_j$ $=$ $e$),
and let  $D$ $=$ $U[x]$.  By (\ref{theo1y}.3),  $P$ $=$ $yU$  is
the maximal ideal in  $U$, so it follows from Proposition \ref{prop1} 
that:
$D$  is a DVR that is a
simple free  integral extension domain of  $U$;
$(QU)D$ $=$ ${P'}^{f_j}$, where
$P'$  $=$ $\pi^{\frac{1}{{f_j}}}D$  is the maximal ideal in  $D$
(so  $D/P'$ $\cong$ $U/P$);
$[D:U]$ $=$ ${f_j}$;
and, $U$  and  $D$  satisfy the Fundamental Equality with no splitting.

Further, since  $u^{\frac{1}{e}}$ $=$ $x$ $\in$ $D$, $D$  and  $V^*$ 
have
the same quotient field  $F(u^{\frac{1}{e}})$.  Thus, $V^*$  and
$D$  are DVRs in  $F(u^{\frac{1}{e}})$, and
$D$ $\le$ $V^*$, so it
follows that  $D$ $=$ $V^*$  and  $P'$ $=$ $N^*$.
Therefore (\ref{theo1y}.4) holds.

It follows from (\ref{theo1y}.3) and (\ref{theo1y}.4) that
$V^*$ $=$ $W[\theta,x]$ ($=$ $W[\theta,u^{\frac{1}{e}}]$)  is a
finite free integral extension domain of  $W$, that  $N^*$
$=$ $xV^*$ $=$ $u^{\frac{1}{e}}V^*$, and that
$[V^* : W]$ $=$ $e_j \cdot f_j$ $=$ $e$.  Also,
since the ramification index of  $U$  over  $W$  is equal to one,
by (\ref{theo1y}.3),
and the ramification index of  $V^*$ $=$ $D$  over  $U$  is equal
to  $f_j$,
by (\ref{theo1y}.4), it follows that
the ramification index of  $V^*$ over  $W$  is equal to  $f_j$.
Further, $[(U/P):(W/Q)]$ $=$ $e_j$, by (\ref{theo1y}.3),
and $V^*/N^*$ $\cong$ $U/P$, by (\ref{theo1y}.4), so
$[(V^*/N^*):(W/Q)]$ $=$ $e_j$.  Thus it follows that
$W$  and  $V^*$  satisfy
the Fundamental Equality with no splitting, hence (\ref{theo1y}.5)
holds.

Moreover, since  $V^*$ is a finite free integral extension
domain of  $W$, and since  $V^*$  is integrally closed, it follows
that  $V^*$  is the integral closure of  $W$ ($=$ $W_j$)  in the 
quotient
field  $F(u^{\frac{1}{e}})$  of  $V^*$.
Therefore it has been shown that if  $V^*$  is an Itoh (e)-valuation 
ring
of  $I$,
then:  $V^*$ $=$ $(\mathbf  T/z^*)_{(p^*/z^*)}$  (for some
(height one) associated prime ideal  $p^*$  of  $u^{\frac{1}{e}} \mathbf 
T$,
where
$z^*$  is the minimal prime ideal in  $\mathbf T$  that is contained
in  $p^*$);   $z$ $=$
$z^* \cap R$  and  $z'$ $=$ $z^* \cap \mathbf R$  are minimal prime 
ideals
in  $R$
and  $\mathbf R$, respectively;
and, if  $F$  (resp., $F(u)$)  is the quotient field
of  $R/z$  (resp., $\mathbf R/z'$), then  $V^*$  is the integral closure
of $W$ ($=$ $V^* \cap F(u)$)  in $F(u^{\frac{1}{e}})$, and  $W$
is a Rees valuation ring of  $u \mathbf R$  (and of  $(u \mathbf R 
+z')/z'$).

It follows that each
Itoh (e)-valuation ring  $V^*$  (with quotient field  
$F(u^{\frac{1}{e}})$)
corresponds to the Rees
valuation ring  $W$ ($=$ $V^* \cap F(u)$)  of  $u \mathbf R$.
On the other hand, if  $W$  is a Rees valuation ring of  $u \mathbf R$,
then:  $W$ $=$ $(\mathbf R'/z')_{(p'/z')}$  for some (height one) 
associated
prime ideal  $p'$  of  $u \mathbf R$, where
$z'$  is the minimal prime ideal in  $\mathbf R$  that is contained
in  $p'$; $z$ $=$
$z' \cap R$  and  $z^*$ $=$
$zR[u^{\frac{1}{e}},t^{\frac{1}{e}}] \cap \mathbf T$
are minimal prime ideals in  $R$
and  $\mathbf T$, respectively;
and, if  $F$  (resp., $F(u^{\frac{1}{e}})$)  is the quotient field
of  $R/z$  (resp., $\mathbf T/z^*$) and  $W''$  is the integral
closure of  $W$  in  $F(u^{\frac{1}{e}})$, then
for each maximal ideal
$M$  in  $W''$, ${W''}_M$  is a Rees valuation ring of
$(u \mathbf T +z^*)/z^*$  and of
$(u^{\frac{1}{e}} \mathbf T +z^*)/z^*$,
by Remarks \ref{RR}.6 and \ref{RR}.5,
so each  ${V''}_M$  is a Rees valuation ring of
$u \mathbf T$ and of  $u^{\frac{1}{e}} \mathbf T $,
by Remark \ref{RR}.3, hence each such ${W''}_M$  is an Itoh 
(e)-valuation
ring of  $I$, by Definition \ref{DD}.
It therefore follows from the first part of this paragraph that  $W''$  
has
exactly one maximal ideal, so the one-to-one correspondence of
(\ref{theo1y}.2) holds.

For (\ref{theo1y}.1), $u^{\frac{1}{e}} V^*$ $=$ $N^*$, by 
(\ref{theo1y}.4).
Therefore, since
$V^*$  is an arbitrary Itoh (e)-valuation ring of  $I$,
it follows from Definitions \ref{DD} and \ref{defi1}.5
that the Rees integers of  $u^{\frac{1}{e}} \mathbf T$
are all equal to one.  Also, since  $\mathbf T$ $=$
$\mathbf R [u^{\frac{1}{e}}]' \cap R[u^{\frac{1}{e}},t^{\frac{1}{e}}]$,
by \cite[Remarks (ii) p. 215]{R2} it follows that  $u^{\frac{1}{e}} 
\mathbf T$ $=$
$(u^{\frac{1}{e}} \mathbf T)_a$.  Therefore it follows
from Remark \ref{RR}.2 that   $\mathbf r$ $=$ $u^{\frac{1}{e}} \mathbf T$
$=$ $(u^{\frac{1}{e}} \mathbf T)_a$ $=$
$\cap \{u^{\frac{1}{e}} {V_j}^* \cap \mathbf T \mid j = 1,\dots,n\}$
$=$ $\cap \{{N_j}^* \cap \mathbf T \mid j = 1,\dots,n\}$  is a 
radical
ideal.
Thus (\ref{theo1y}.1) holds, so
(\ref{theo1y}.1) - (\ref{theo1y}.5) hold.

Finally, for (\ref{theo1y}.6), let  $W$ be a
Rees valuation ring of  $u \mathbf R$,
let  $V^*$  be the corresponding (as in (\ref{theo1y}.2)) Itoh 
(e)-valuation
ring of  $I$,
and let  $U$  be as
in (\ref{theo1y}.3), so  $W$ $\le$ $U$ $\le$
$V^*$.  Then there exists  $j$ $\in$ $\{1,\dots,n\}$  such that
$[U:W]$ $=$ $e_j$  (by (\ref{theo1y}.3)) $=$ $e$ (by the hypothesis
in (\ref{theo1y}.6)) $=$  $[V^*:W]$ (by (\ref{theo1y}.5)), so
$V^*$ $=$ $U$  is a simple free integral extension domain of  $W$.
Also, $P$ $=$ $N^*$ $=$ $u^{\frac{1}{e}}V^*$ (by (\ref{theo1y}.5)) 
$= $
$u^{\frac{1}{e}}U$.
Since this holds for all Itoh (e)-valuation rings $(U,P)$
of  $I$,  $u^{\frac{1}{e}} \mathbf T$
is a radical ideal, hence (\ref{theo1y}.6) holds.
\end{proof}

Concerning Theorem \ref{theo1y}.6, it is shown in Corollary \ref{co2} 
below
that
there always exists a finite integral extension ring  $A$  of  $R$  such
that the Rees integers of  $IA$  are all equal to  $e$, where
$e$  is an arbitrary common multiple of the Rees integers
$e_1,\dots,e_n$  of  $I$, so Theorem \ref{theo1y}.6 directly
applies to  $IA$  in place of  $I$.

\begin{rema}
\label{x1}
{\em
Let  $(V,N)$  be a Rees valuation ring of  $I$  and assume
that $IV$ $=$ $N^k$.  Let  $(W,Q)$
be the Rees valuation ring of  $uR[u,tI]$  that
corresponds (as in Proposition \ref{prop3}) to  $V$, and
let  $h$  be an arbitrary positive integer.
Then it follows from Theorem \ref{theo1y}.3 and
\ref{theo1y}.4 that:
$U$ $=$ $W[u^{\frac{1}{k}}]'$  and
$D$ $=$ $W[u^{\frac{1}{hk}}]'$ $=$ 
$U[(u^{\frac{1}{k}})^{\frac{1}{h}}]' $
are DVRs that are finite free integral extension domains
of  $W$;
the maximal ideal of  $U$  (resp., $D$)  is
$P$ $=$ $u^{\frac{1}{k}} U$
(resp.,  $M$ $=$ $u^{\frac{1}{hk}} D$);
the ramification index of  $D$  relative to  $U$  (resp.,
$U$  relative to  $W$)  is equal to  $h$  (resp., $1$);
$[D:U]$ $=$ $h$ and  $[U:W]$ $=$ $k$; and,
$D/M$ $\cong$ $U/P$ and $[(U/P):(W/Q)]$ $=$ $k$.
It therefore follows that  $W$ and  $U$  (resp., $W$ and  $D$,
$U$  and  $D$) satisfy the
Fundamental Equality with no splitting.
Also, it follows from Theorem \ref{theo1y}.3 (resp., \ref{theo1y}.4)
that  $U$  (resp., $D$)  is an Itoh (k)-valuation (resp.,
(hk)-valuation) ring of  $I$.
}

\end{rema}

\begin{rema}
\label{Independent}
{\em
The following result is called the Theorem of
Independence of Valuations, and it is proved in
\cite[(11.11)]{N2}:  Let  $(V_1,N_1),\dots,(V_n,N_n)$
be valuation rings with the same quotient field  $F$, and
assume there are no containment relations among the  $V_j$.
Let  $\mathbf V$ $=$ $V_1 \cap \cdots \cap V_n$  and let
$P_j$ $=$ $N_j \cap \mathbf V$  ($j$ $=$ $1,\dots,n$).  Then
$P_1,\dots,P_n$  are the maximal ideals in  $\mathbf V$
and  $V_j$ $=$ $\mathbf V_{P_j}$  for  $j$ $=$ $1,\dots,n$.
}

\end{rema}

\begin{coro}
\label{co1y}
Let  $I$  $=$ $(b_1,\dots,b_g)R$  be a nonzero proper ideal in a
Noetherian integral domain
$R$, let  $F$  be the quotient field of  $R$, let  $\mathbf R$ $=$
$R[u,tI]$, let  $(W_1,Q_1),\dots,(W_n,Q_n)$  be the Rees valuation
rings of  $u \mathbf R$, let  $e$  be a positive common multiple
of the Rees integers  $e_1,\dots,e_n$   of  $u \mathbf R$, let
$\mathbf W$ $=$ $W_1 \cap \cdots \cap W_n$, and let
$\mathbf D$ $=$ $\mathbf W[x]$, where  $x$ $=$ $u^{\frac{1}{e}}$
in a fixed algebraic closure  $F(u)^*$  of  $F(u)$.  Then:
\begin{enumerate}[$(1)$]

\item
$\mathbf W$  is a semi-local Dedekind domain with exactly
$n$  maximal ideals $P_j$ $=$ $Q_j \cap \mathbf W$, $j$ $=$ 
$1,\dots,n$.

\item
$\mathbf D$ is a simple free integral extension domain of
$\mathbf W$  of rank  $e$, and  $\mathbf D$  has exactly  $n$  maximal
ideals
$M_j$ $=$ $(P_j,x) \mathbf D$.

\item
There exist distinct elements  $\theta_1,\dots,\theta_n$
in the quotient field  $F(x)$  of  $\mathbf D$  such that:

%\newline
\hspace*{.3in}{\bf (3.1)}
$\mathbf D'$  is the intersection of the  Itoh 
\textup{(}e\textup{)}-valuation
rings  $({V_1}^*,{N_1}^*),\dots,({V_n}^*,{N_n}^*)$
\hspace*{.6in}  of  $I$,
where  ${V_j}^*$ $=$ $W_j[\theta_j,x]$ $=$
${\mathbf D'}_{\mathbb {M}_j }$, where
$\mathbb {M}_j$ $=$  ${N_j}^* \cap \mathbf D'$.

%\newline
\hspace*{.3in}{\bf (3.2)}
$\mathbf D'$ $=$ $\mathbf D[ \theta_1,\dots,\theta_n]$  is a
semi-local Dedekind domain that is a
finite integral
\hspace*{.6in} extension domain of  $\mathbf D$,
and  $\mathbf D'$  has exactly  $n$  maximal ideals
$\mathbb {M}_j$ $=$  ${N_j}^* \cap \mathbf D'$
\hspace*{.6in} 
$(P_j,x,{\theta_1},$ 
$\dots,{\theta_{j-1}},{\theta_{j+1}},\dots,{\theta_n})
\mathbf D'$.

%\newline
\hspace*{.3in}
{\bf (3.3)}
The Jacobson radical of  $\mathbf D'$  is  $x \mathbf D'$.

\item
Assume that  $e_1$ $=$ $\cdots$ $=$ $e_n$ $=$ $e$  and that there 
exists
$b$ $\in$ $I$  such that  $bW_j$ $=$ $IW_j$  for  $j$
$=$ $1,\dots,n$.  Then, with  $\mathbf D'$, as in 
\textup{(}3.1\textup{)},
and $Q_j$ $=$ $\pi_j W_j$  for  $j$ $=$ $1,\dots,n$,
$\frac{x}{\pi_1 \cdots \pi_n}$  is a unit in  $\mathbf D'$,
$\mathbf D'$ $=$ $\mathbf D[\frac{x}{\pi_1 \cdots \pi_n}]$
$=$ $\mathbf W[\frac{x}{\pi_1 \cdots \pi_n}]$,
and  $x \mathbf D'$  is
the Jacobson radical of  $\mathbf D'$.
\end{enumerate}

\end{coro}

\begin{proof}
(1) and part of (3.1) follow immediately from the Independence
of Valuations Theorem (see Remark \ref{Independent}).
The first part of (2) is clear.  Also,
each  $M_j$  is a maximal
ideal in  $\mathbf D$, since  $\mathbf D/M_j$ $\cong$ $W_j/Q_j$.
Further, $u$  (resp., $x$)
is in the Jacobson radical of  $\mathbf W$  (resp., $\mathbf D$),
$x \mathbf D \cap \mathbf W$ $=$ $u \mathbf W$  (since  $\mathbf D$  
is a
free $\mathbf W$-algebra),
and  $\mathbf D/(x \mathbf D)$ $\cong$ $\mathbf W/(u \mathbf W)$,
so it follows that   $\mathbf W$,
$\mathbf W /(u \mathbf W)$,
$\mathbf D$, and $\mathbf D /(x \mathbf D)$
each have exactly  $n$  maximal ideals,
so (2) holds.

For (3), for  $j = 1,\dots,n$,  let  $(V_j,N_j)$ be
the Rees valuation ring of  $I$  that corresponds
(as in Proposition \ref{prop3}) to  $(W_j,Q_j)$
(so there exists
$b_{\sigma(j)}$ (where  $\sigma~:~\{1,\dots,n\} 
~\rightarrow~\{1,\dots,g\}$)
in  $\{b_1,\dots,b_g\}$  such that
$V_j$ $=$ $(R[\frac{I}{b_{\sigma(j)}}]')_{p_j}$, for some
(height one) prime ideal  $p_j$  in $R[\frac{I}{b_{\sigma(j)}}]'$),
so  $W_j$ $=$ $V_j[tb_{\sigma(j)}]_{N_jV_j[tb_{\sigma(j)}]}$,
$W_j \cap F$ $=$ $V_j$, and  $Q_j \cap F$ $=$ $N_j$.

Let  $\mathbf V$ $=$ $V_1 \cap \cdots \cap V_n$, and
for  $j = 1\dots,n$,  let $q_j$ $=$ $N_j \cap \mathbf V$.
Then   $\mathbf V$ $=$ $\mathbf W \cap F$  and  $\mathbf W$ $=$
$\mathbf V[u]_S$, where  $S$ $=$
$\mathbf V[u] - (q_1 \mathbf V[u] \cup \cdots \cup q_n \mathbf V[u])$,
and  $q_j$ $=$ $(Q_j \cap F) \cap \mathbf V$ $=$ $(Q_j \cap \mathbf 
W)
\cap \mathbf V$
$=$ $P_j \cap \mathbf V$, for  $j$ $=$ $1,\dots,n$.
Therefore it follows from the Independence of Valuations Theorem
(see Remark \ref{Independent})
that $\mathbf V$  is a semi-local Dedekind domain with exactly
$n$  maximal ideals  $q_j$ ($j$ $=$ $1,\dots,n$).
Thus  $\mathbf V$  is a Principal Ideal Domain,
by \cite[Theorem 16, p. 278 ]{ZS1},
so for  $j = 1,\dots,n$,  there exists  $\pi_j$ $\in$ $q_j$
such that  $q_j$ $=$ $\pi_j \mathbf V$, so:
$\pi_j V_j$ $=$ $(N_j \cap \mathbf V) \mathbf V_{q_j}$ $=$ $N_j$;
$\pi_j \mathbf W$ $=$ $q_j \mathbf W$ $=$ $(P_j \cap \mathbf V) 
\mathbf W$
$=$ $P_j$; and,
$\pi_j W_j$ $=$ $q_j \mathbf W_{Q_j \cap \mathbf W}$
$=$ $P_j W_j$ $=$ $Q_j$.
Then, since the Rees integers of  $u \mathbf R$  are  $e_1,\dots,e_n$
(by hypothesis),
$u \mathbf W$ $=$ ${P_1}^{e_1} \cap \cdots \cap {P_n}^{e_n}$ $=$ 
${P_1}^{e_1} \cdots {P_n}^{e_n}$ $=$
${\pi_1}^{e_1} \cdots {\pi_n}^{e_n} \mathbf W$, so
there exists a unit  $w$ $\in$ $\mathbf W$  such that
$u$ $=$ $w{\pi_1}^{e_1} \cdots {\pi_n}^{e_n}$,
hence
$$(\textup{*}1)~~~~~~~~~~\text{there exists a unit}~w_j \in 
W_j~~\text{such that}~u
= w_j{\pi_j}^{e_j}~~\text{(in}~  W_j = \mathbf W_{P_j}\text{), for} 
~j
= 1,\dots,n,$$
where
$$(\textup{*}2)~~~~~~~w_j = \frac{u}{\pi_j^{e_j}} = w \prod_{i \neq 
j} \pi_i^{e_i}
\in (\bigcap_{i \neq j} P_i) - P_j,~\text{for }~j = 1,\dots,n.$$

Also, since  $e_j$  is the Rees integer of  $I$  with respect to  $V_j$  
and
$IV_j$ $=$ $b_{\sigma(j)}V_j$, there exists
a unit  $v_j$ $\in$ $V_j$ such that
$$(\textup{*}3)~~~~~~~~~~b_{\sigma(j)} = v_j {\pi_j}^{e_j} ~~ 
\text{in}~~
V_j\text{, for}~~  j = 1,\dots,n.$$

Therefore, since  $b_{\sigma(j)}$ $=$ $u(tb_{\sigma(j)})$  (in  
$\mathbf R$
$\subseteq$ $W_j$), it follows from (*3) and (*1)
that
$v_j {\pi_j}^{e_j}$$=$ $b_{\sigma(j)}$ $=$ $u(tb_{\sigma(j)})$ $=$ 
$w
_j {\pi_j}^{e_j} (tb_{\sigma(j)})$,
hence
\begin{equation*}
(\textup{*}4)\quad w_j~ = \frac{v_j}{tb_{\sigma(j)}} \text{ and } w_j 
+ P_j  \text{ is transcendental over }\mathbf V/(P_j \cap \mathbf V)  
\text{ for } j = 1,\dots,n.
\end{equation*}
Therefore
$$(\textup{*}5)\quad \overline m_{e_j}(X) = X^{e_j} - (w_j+P_j) \text{ 
is irreducible in } (\mathbf W/P_j)[X]~~  for~~  j = 1,\dots,n.$$

For  $j$ $=$ $1,\dots,n$,  let  $\theta_j$ $=$ 
${w_j}^{\frac{1}{e_j}}$
in the (fixed) algebraic closure  $F(u)^*$  of  $F(u)$.
Then it is shown
in Theorem \ref{theo1y}.3 - \ref{theo1y}.5 that, for
the Itoh (e)-valuation ring  ${V_j}^*$  of  $I$,
$$(\textup{*}6)\quad {V_j}^* = W_j[\theta_j,x] \mbox{ and } \theta_j
\mbox{ is a unit in }
{V_j}^*, \mbox{ for } j = 1,\dots,n.$$
Also, ${V_1}^* \cap \cdots \cap {V_n}^*$ is the integral
closure $\mathbf W^*$  of  $\mathbf W$  in  $F(x)$,
by Theorem \ref{theo1y}.2, and  $x$ $=$
$u^{\frac{1}{e}}$ $\in$ $\mathbf W^*$  (since  $u$
$\in$ $\mathbf W$), so it follows
that  $\mathbf W^*$ $=$ $\mathbf W[x]'$ $=$ $\mathbf D'$.

Further, for  $j$ $=$ $1,\dots,n$,
$\theta_j$ $=$ ${w_j}^{\frac{1}{e_j}}$ $\in$ $\mathbf W^*$,
by (*2) and (*6), and by (*2)
$$\theta_1,\dots,\theta_{j-1},\theta_{j+1},\dots,\theta_n \in {N_j}^*
\cap \mathbf W^*
= {N_j}^* \cap \mathbf D' = \mathbb M_j,
\mbox{ for }  j = 1,\ldots,n.$$

Hence by (*2) and (*6) we see that
$$(\textup{*}7)\quad \text{for} ~~i~ = 1,\dots,n, ~{\theta_i} \in
\cap \{{N_j}^* \cap \mathbf D' \mid j = 1,\dots,i-1,i+1,\dots,n\} - 
({N_i}^*
\cap \mathbf D').$$

Let  $\mathbf E$ $=$ $\mathbf D[ {\theta_1},\dots,{\theta_n}]$. Then,
for  $j$ $=$ $1,\dots,n$,
$x$ $\in$ $\mathbf D$ $\subseteq$ $\mathbf E$
$\subseteq$ $\mathbf W^*$ (by (*2) and (*6))
$\subseteq$ ${V_j}^*$
$=$ $W_j[\theta_j ,x]$.
Therefore  $W_j$ $=$ ${\mathbf W}_{{N_j}^* \cap \mathbf W}$ 
$\subseteq$
$\mathbf D_{{N_j}^* \cap \mathbf D}$ $\subseteq$
$\mathbf E_{{N_j}^* \cap \mathbf E}$ $\subseteq$
${\mathbf W^*}_{{N_j}^* \cap \mathbf W^*}$
$=$ ${V_j}^*$
$=$ $(W_j[x])[\theta_j]$ $\subseteq$
$\mathbf E_{{N_j}^* \cap \mathbf E}$, so
${V_j}^*$ $=$ $\mathbf E_{{N_j}^* \cap \mathbf E}$, for
$j$ $=$ $1,\dots,n$. Also,
$\mathbf D$ and  $\mathbf W^*$  have exactly  $n$  maximal ideals,
so it follows from integral dependence that  $\mathbf E$  has exactly
$n$ maximal ideals.  Further, for
each integral domain  $A$, $A$ $=$ $\cap \{A_M \mid M$  is a maximal
ideal in  $A\}$, by \cite[(33.9)]{N2}.  Therefore  it follows that
$\mathbf E$ $=$ $\mathbf D[{\theta_1},\dots,{\theta_n}]$ $=$
${V_1}^* \cap \cdots \cap {V_n}^*$ $=$ $\mathbf D'$.
Thus (3.1) and (3.2) hold.

(3.3) follows immediately from Theorem \ref{theo1y}.4.

Finally, for (4),
let  $(V_j,N_j)$  be
the Rees valuation ring of  $I$  that corresponds
(as in Proposition \ref{prop3}) to  $(W_j,Q_j)$.
Then it follows from the hypothesis on  $b$
that, for  $j$ $=$ $1,\dots,n$,  $W_j$ $=$ $V_j[tb]_{N_jV_j[tb]}$  
and
$Q_j$ $=$ $N_jW_j$.

Let  $\mathbf V$ $=$ $V_1 \cap \cdots \cap V_n$,
so by Remark \ref{Independent}  $\mathbf V$  is a semi-local Dedekind
domain with exactly  $n$  maximal ideals
$N_j \cap \mathbf V$ ($=$ $P_j \cap \mathbf V$).
Thus  $\mathbf V$  is a Principal Ideal Domain,
by \cite[Theorem 16, p. 278]{ZS1},
so for  $j$ $=$ $1,\dots,n$  there exists  $\pi_j$ $\in$ $N_j \cap 
\mathbf V$
such that  $\pi_j \mathbf V$ $=$ $N_j \cap \mathbf V$.

Let  $\alpha$ $=$ $\pi_1 \cdots \pi_n$,
let  $\mathbf N$ $=$ $N_1 \cap \cdots \cap N_n$,
and let  $\mathbf Q$ $=$ $Q_1 \cap \cdots \cap Q_n$.  Then,
since the Rees integers of  $I$  are all equal to  $e$,
it follows that
$$(\ref{co1y}.4.1) \; \; \; I \mathbf V = b \mathbf V = \mathbf N ^e
= \alpha^e \mathbf V \mbox{ and } I \mathbf W = u \mathbf W = \mathbf 
Q ^e
= \alpha^e \mathbf W.$$

It follows from (\ref{co1y}.4.1) that
$$(\ref{co1y}.4.2) \mbox{ there exist units } v \in \mathbf V \mbox{ and } w 
\in \mathbf W \mbox{ such that } b = v \alpha^e  \in  \mathbf V \mbox{ and }
u = w \alpha^e  \in  \mathbf W.$$

Since
$b$ $=$ $u(tb)$  in  $\mathbf R$ $\subseteq$ $\mathbf W$,
and since  $Q_j = N_jV_j[tb]_{N_jV_j[tb]}$  for 
$j = 1, \ldots, n$, it follows from (\ref{co1y}.4.2)
that
$$(\ref{co1y}.4.3)~~  w = \frac{v}{tb} \mbox{ and } w + P_j 
\mbox{ is transcendental over } \mathbf V/(P_j \cap \mathbf V)
\mbox{ for } j = 1,\dots,n.$$

Therefore
$$(\ref{co1y}.4.4)~~ \overline m_e(X) = X^e - (w+P_j) \mbox{ is 
irreducible in } (\mathbf W/P_j)[X] \mbox{ for }  j = 1,\dots,n.$$

Let  $\theta$ $=$ $w^{\frac{1}{e}}$  in the
(fixed) algebraic closure  $F(u)^*$  of  $F(u)$.  Then
it follows from Proposition \ref{prop2}, together with
(\ref{co1y}.4.4),
that:
$m_e(X)$  is irreducible in the UFD  $\mathbf W[X]$;
for  $j$ $=$ $1,\dots,n$, $\mathbb M_j$ $=$ $P_j \mathbf E$ is a
(principal) maximal ideal,
where  $\mathbf E$ $=$ $\mathbf W[\theta]$;
and, $[\mathbf E: \mathbf W]$ $=$ $e$.

Also, $\mathbf E$  is integral over  $\mathbf W$, so it follows
that  $\mathbb M_1,\dots,\mathbb M_n$  are the only nonzero prime ideals 
in
$\mathbf E$,
hence  $\mathbf E$ is a semi-local Dedekind domain that is a simple free
integral extension
domain of  $\mathbf W$  and $\mathbf Q \mathbf E$ $=$
$(Q_1 \cap \cdots \cap Q_n)\mathbf E$ $=$
$(\mathbb M_1 \cap \cdots \cap \mathbb M_n)\mathbf E$ $=$ $\pi_1 
\cdots
\pi_n \mathbf E$ $=$
$\alpha \mathbf E$  is its Jacobson radical  $J$.  Therefore  $\mathbf 
D'$
$=$
$\mathbf E$ $=$ $\mathbf W[\theta]$ $=$ $\mathbf D[\theta]$, and 
since
$w$ $=$ $\frac{u}{\alpha ^e}$  is a unit in  $\mathbf W$, it follows
that $\theta$ $=$ $w^{\frac{1}{e}}$
$=$ $\frac{x}{\alpha}$  is a unit in  $\mathbf D'$  and that
$J$ $=$ $\alpha \mathbf E$ $=$ $x\mathbf E$.
\end{proof}

The next remark lists several well known
facts concerning finite field extensions and
ramification.

\begin{rema}
\label{remxx}
{\em
Let  $(U_1,P_1)$ $\le$
$(U_2,P_2)$ $\le$
$(U_3,P_3)$ be DVRs such that
$U_3$  is a finite integral extension of  $U_1$.
Then:

\noindent
{\bf (1)} $[(U_3)_{(0)} : (U_1)_{(0)}]$ $=$ $[(U_3)_{(0)} : 
(U_2)_{(0)}]$ $\cdot$ $
[(U_2)_{(0)} : (U_1)_{(0)}]$.

\noindent
{\bf (2)} $[(U_3/P_3):(U_1/P_1)]$ $=$
$[(U_3/P_3):(U_2/P_2)]$ $\cdot$ $[(U_2/P_2):(U_1/P_1)]$.

\noindent
Also, if we let  $r_{3,1}$  (resp., $r_{3,2}$, $r_{2,1}$)
denote the ramification index of  $U_3$ relative to  $U_1$
(resp., $U_3$  relative to  $U_2$, $U_2$ relative to  $U_1$),
then:

\noindent
{\bf (3)} $r_{3,1}$ $=$ $r_{3,2} \cdot r_{2,1}$.

\noindent
Further, if  $U_1$ and  $U_3$
satisfy the Fundamental Equality with no splitting, then

\noindent
{\bf (4)} $r_{3,1} \cdot [(U_3/P_3):(U_1/P_1)]$ $=$ $[(U_3)_{(0)} : 
(U_1)_{(0)}]$.

\noindent
It then follows from (1) - (4)
and the Fundamental Inequality, (\ref{FI}), that:

\noindent
{\bf (5)} $r_{3,2} \cdot [(U_3/P_3):(U_2/P_2)]$ $=$ $[(U_3)_{(0)} : 
(U_2)_{(0)}]$.

\noindent
{\bf (6)} $r_{2,1} \cdot [(U_2/P_2):(U_1/P_1)]$ $=$ $[(U_2)_{(0)} : 
(U_1)_{(0)}]$.

\noindent
{\bf (7)} Hence, if $U_1$ and  $U_3$
satisfy the Fundamental Equality with no splitting, then
both  ($U_1$ and  $U_2$) and ($U_2$  and  $U_3$)
satisfy the Fundamental Equality with no splitting.
}

\end{rema}

The next proposition shows that Theorem \ref{theo1y}.2 holds for
all integers greater than one.

\begin{prop}
\label{ptheo1y}
With the notation of Theorem \ref{theo1y}, let  $k$ $\ge$ $2$  be
an arbitrary integer, let  $\mathbf S_k$
$=$ $\mathbf R[u^{\frac{1}{k}}]$, and let
$\mathbf T_k$ $=$ ${\mathbf S_k}' \cap 
R[u^{\frac{1}{k}},t^{\frac{1}{k}}]$.
Then:

\begin{enumerate}[$(1)$]
\item There exists a one-to-one correspondence between
the Itoh \textup{(}k\textup{)}-valuation rings  $(U,P)$  of  $I$
and the Rees valuation rings  $(W,Q)$  of
$u \mathbf R$; namely, given  $W$, if  $F(u)$  is the quotient
field of  $W$, then  $U$  is the integral closure
of  $W$  in  $F(u^{\frac{1}{k}})$.

\item Let $W$ and $U$  be corresponding \textup{(}as in
\textup{(}1\textup{))}.  Then
$U$   is a finite integral extension domain of  $W$, and
$W$  and  $U$  satisfy the Fundamental Equality
with no splitting,
\end{enumerate}

\end{prop}

\begin{proof}
Let  $m$  be the least common multiple of
$e_1,\dots,e_n$, and let  $e$ $=$ $k \cdot m$.
Then it is clear that
$\mathbf R$ $\subseteq$
$\mathbf T_k$ $\subseteq$
$\mathbf T_e$ $=$
$\mathbf R[u^{\frac{1}{e}}]' \cap R[u^{\frac{1}{e}},t^{\frac{1}{e}}]$
and that $\mathbf T_e$  is integral over  $\mathbf R$.
Also, by Theorem \ref{theo1y}.2
there exists a one-to-one correspondence between
the Rees valuation rings  $(W,Q)$  of  $u \mathbf R$  and the
Itoh (e)-valuation rings  $(V^*,N^*)$  of  $I$,
and  $V^*$  is the integral closure of  $W$  in the quotient field
of  $V^*$. It follows from this, and
integral dependence, that (1) holds.

For (2), it follows from the proof of (1) that if  $W$  and
$V^*$  are as in (1), then there exists an Itoh (k)-valuation
ring  $U$  of  $I$  such that  $W$ $\le$ $U$ $\le$ $V^*$.
Also, $V^*$  is a finite free integral extension domain of  $W$,
by Theorem \ref{theo1y}.5, hence  $U$
is a finite integral extension domain of  $W$.
Further, $W$  and  $V^*$  satisfy the Fundamental Equality with no
splitting, by Theorem \ref{theo1y}.5, so
it follows from Remark \ref{remxx}.7 that (2) holds.
\end{proof}

\begin{term}
\label{t2}
{\em
If  $(W,Q)$  is a Rees valuation ring of
$uR[u,tI]$  and  $(U,P)$  is the
corresponding (as in
Proposition \ref{ptheo1y}.1)
Itoh (k)-valuation ring of  $I$
(so  $U$ $\ge$ $W$  and  $U$  is the
integral closure of  $W$  in the quotient
field  $F(u^{\frac{1}{k}})$  of  $U$),
then we say that
$\mathbf{(U,P)}$
{\bf{is the Itoh (k)-valuation ring of}}
$I$
{\bf{that corresponds to}}
$\mathbf{(W,Q)}$.  Also, if  $(V,N)$  is
the Rees valuation ring of  $I$  that corresponds
(as in Proposition \ref{prop3}) to  $(W,Q)$, then
we say that $(\mathbf{U,P)}$
{\bf{is the Itoh (k)-valuation ring of}}
$I$
{\bf{that corresponds to}}
$(\mathbf{V,N)}$.
}

\end{term}

\begin{prop}
\label{prop6}
Let  $I$  be a regular proper ideal in a Noetherian ring  $R$,
let  $(V,N)$  be a Rees valuation ring of  $I$, let  $(W,Q)$  be
the Rees valuation ring of  $uR[u,tI]$  that corresponds
\textup{(}as in Proposition \ref{prop3}\textup{)} to  $(V,N)$,
and let  $(U_k,P_k)$  be the Itoh \textup{(}k\textup{)}-valuation ring 
of
$I$  that corresponds
to  $(V,N)$  \textup{(}see \textup{(}\ref{t2}\textup{))}, where  $k$ 
$\ge$
$2$  is an arbitrary
integer.  Assume that  $IV$ $=$ $N^e$.  Then the following hold:

\begin{enumerate}[$(1)$]
\item If  $e$  is a multiple of  $k$, then  $P_k$ $=$ $NU_k$  and
$[(U_k)_{(0)} : W_{(0)}]$ $=$ $k$ $=$ $[(U_k/P_k):(W/Q)]$.

\item If  $e$  and  $k$  are relatively prime, then  $NU_k$ $=$  
${P_k}^k$,
$[(U_k)_{(0)} : W_{(0)}]$ $=$ $k$, and  $U_k/P_k$ $\cong$ $W/Q$.

\item If  the greatest common divisor of  $e$  and  $k$  is  $d$,
and if  $c$ $\in$ $\mathbb N_{> 0}$  is such that $cd$ $=$ $k$,
then  $NU_k$ $=$  ${P_k}^c$,
$[(U_k)_{(0)} : W_{(0)}]$ $=$ $k$,
and  $[(U_k/P_k): (W/Q)]$ $=$ $d$.
\end{enumerate}

\end{prop}

\begin{proof}
For (1), let  $(U_e,P_e)$  be the Itoh (e)-valuation ring of  $I$
that corresponds to  $(V,N)$.  Then it follows from Remark \ref{x1}
that  $P_e$ $=$ $NU_e$  and
$[(U_e)_{(0)} : W_{(0)}]$ $=$ $e$ $=$ $[(U_e/P_e):(W/Q)]$.  Also, it 
is
clear
that  $W$ $\le$ $U_k$ $\le$ $U_e$  and that
$[(U_k)_{(0)}:W_{(0)}]$ $=$ $k$. Item (1) readily follows from this, 
together
with Remark \ref{remxx}.1 - \ref{remxx}.3.

For (2), let  $(U_{ke},P_{ke})$  be the Itoh (ke)-valuation
ring of  $I$  that corresponds to  $(V,N)$.
Then it follows from Remark \ref{x1}
that  ${P_{ke}}^k$ $=$ $NU_{ke}$,
$[(U_{ke})_{(0)} : W_{(0)}]$ $=$ $ke$, and
$[(U_{ke}/P_{ke}):(W/Q)]$ $=$ $e$.  Also, it is clear
that  $W$ $\le$ $U_k$ $\le$ $U_{ke}$  and that
$[(U_k)_{(0)}:W_{(0)}]$ $=$ $k$.  Since  $e$  and  $k$  are relatively
prime,
(2) readily follows from this, together
with Remark \ref{remxx}.2, \ref{remxx}.3, and \ref{remxx}.7.
\newpage

For (3), let  $(U_c,P_c)$  (resp., $(U_d,P_d)$)
be the Itoh (c)-valuation  (resp., (d)-valuation)
ring of  $I$  that corresponds to  $(V,N)$.
Since  $cd$ $=$ $k$, it follows that  $W$ $\le$ $U_c$ $\le$ $U_k$
and  $W$ $\le$ $U_d$ $\le$ $U_k$.  Since  $e$  is a multiple of  $d$,
it follows from (1) that
$$(a)~~~~~~~~~~ QU_d = P_d~~  \text{and} ~~ [(U_d)_{(0)} : W_{(0)}] 
= d
= [(U_d/P_d):(W/Q)].$$
Since  $c$  and  $e$  are
relatively prime, it follows from (2) that
$$(b)~~~~~~~~~~ QU_c = {P_c}^c~~ and ~~U_c/P_c \cong W/Q.$$
It follows from (a), (b), and Remark \ref{remxx}.2, \ref{remxx}.3,
and \ref{remxx}.7 that (3) holds.
\end{proof}

Proposition \ref{prop5} gives several equivalences 
to Property
\ref{theo1y}.1
of Theorem \ref{theo1y}
and also shows that
the hypothesis
of Theorem \ref{theo1y}
is necessary for
\ref{theo1y}.1.

\begin{prop}
\label{prop5}
Let  $I$  be a regular proper ideal in a Noetherian ring  $R$,
let  $\mathbf R$ $=$ $R[u,tI]$  be the
Rees ring of  $R$  with respect to  $I$,
let  $k$ $\ge$ $2$  be an arbitrary integer,
let  $\mathbf S_k$
$=$ $\mathbf R[u^{\frac{1}{k}}]$, and let  $\mathbf T_k$
$=$ ${\mathbf {S_k}}' \cap R[u^{\frac{1}{k}},t^{\frac{1}{k}}]$.
Then the following statements are equivalent:

\begin{enumerate}[$(1)$]
\item $u^{\frac{1}{k}} \mathbf T_k$  is a radical ideal.

\item
The Rees integers of  $u^{\frac{1}{k}} \mathbf T_k$  are all equal to 
one.

\item
The Rees integers of  $(u^{\frac{1}{k}} \mathbf S_k)_a$  are all equal 
to
one.

\item
$k$  is a common multiple of the Rees integers of  $I$.
\end{enumerate}

\end{prop}

\begin{proof}
A primary decomposition of the regular proper principal radical ideal
$u^{\frac{1}{k}} \mathbf T_k$, together with
Remark \ref{RR}.4, shows that (1) $\Rightarrow$ (2),
and the next
to last paragraph of the proof of Theorem \ref{theo1y}
shows that (2) $\Rightarrow$ (1).

Also, (2) $\Leftrightarrow$ (3),
by Remark \ref{RR}.5, since  $\mathbf S_k$ $\subseteq$ $\mathbf T_k$
$\subseteq$ ${\mathbf S_k}'$.

Assume that (2) holds and let  
$({V_1}^*,{N_1}^*),\dots,({V_n}^*,{N_n}^*)$  be
the Itoh (k)-valuation rings of  $I$, so  $u^{\frac{1}{k}}{V_j}^*$ $=$
${ N_j}^*$  for
$j$ $=$ $1,\dots,n$.  For  $j$ $=$ $1,\dots,n$,  let  $uW_j$ $=$ 
${Q_j}^{e_j}$,
where  $\mathbf {RV}(uR[u,tI])$ $=$ $\{(W_j,Q_j) \mid j = 
1,\dots,n\}$.
Suppose that  $k$  is not a multiple of  $e_j$  for some  $j$, let
$d$  be the greatest common divisor of  $k$  and  $e_j$, and let
$c \ge 1$  and  $h > 1$  be integers  such that  $cd$ $=$ $k$  and
$hd$ $=$ $e_j$.  Then it follows from Proposition \ref{prop6}(3)
(with  $(W_j,Q_j)$  (resp., $({V_j}^*,{N_j}^*)$)  in place of
$(W,Q)$  (resp., $(U_k,P_k)$)  that
$Q_j {V_j}^*$ $=$ $({N_j}^*)^c$.  Since  $uW_j$ $=$ ${Q_j}^{e_j}$,
it follows that  $u^{\frac{1}{k}}{V_j}^*$ $=$ 
$({N_j}^*)^{\frac{ce_j}{k}}$
$=$ $({N_j}^*)^{h}$.  However, $h$ $>$ $1$, and this
contradicts (2).  Therefore the supposition that  $k$ is not a multiple
of  $e_j$  leads to a contradiction, hence (2) $\Rightarrow$ (4).

Finally, (4) $\Rightarrow$ (1) by
Theorem
\ref{theo1y}.1.
\end{proof}

\begin{coro}
\label{end}
Let  $I$  be a regular proper ideal in a Noetherian ring  $R$,
let  $\mathbf R$ $=$ $R[u,tI]$  be the
Rees ring of  $R$  with respect to  $I$,
and for each integer  $k$ $\ge$ $2$
let  $\mathbf S_k$
$=$ $\mathbf R[u^{\frac{1}{k}}]$ and let  $\mathbf T_k$
$=$ ${\mathbf {S_k}}' \cap R[u^{\frac{1}{k}},t^{\frac{1}{k}}]$.
Then the following statements are equivalent:

\noindent
{\bf (1)}
The Rees integers of  $I$  are all equal to one.

\noindent
{\bf (2)}
For all integers  $k$ $\ge$ $2$, the ideal  $u^{\frac{1}{k}} \mathbf T_k$
is a radical ideal.

\end{coro}

\begin{proof}
This follows immediately from Proposition \ref{prop5}(1) 
$\Leftrightarrow$
(4).
\end{proof}

\section{A RELATED THEOREM}
In this section, we first prove an expanded version of
Corollary \ref{co1y} and Theorem \ref{theo1y}.6,
and we then prove a closely related and more general theorem.

The next theorem is an expanded version of
Corollary \ref{co1y} and Theorem \ref{theo1y}.6.
One of the main reasons for including  this theorem is to
show that it displays, except for separability,
a realization (see Definition \ref{real}
below) of a powerful classical method.  This is explained somewhat
more fully in Remark \ref{EXP} below.

\begin{theo}
\label{theox}
Let  $I$  be a nonzero proper ideal in a Noetherian integral domain
$R$, let  $F$  be the quotient field of  $R$, let  $\mathbf R$ $=$
$R[u,tI]$, let  $(W_1,Q_1),\dots,(W_n,Q_n)$  be the Rees valuation
rings of  $u \mathbf R$, let  $e$  be a positive common multiple
of the Rees integers  $e_1,\dots,e_n$   of  $u \mathbf R$, let
$\mathbf W$ $=$ $W_1 \cap \cdots \cap W_n$, let  $\mathbf Q$ $=$
$Q_1 \cap \dots \cap Q_n$, and for  $j$ $=$ $1,\dots,n$,  let  $P_j$
$=$ $Q_j \cap \mathbf W$, so  $\mathbf W$  is a semi-local Dedekind
domain, $\mathbf Q$ is its Jacobson radical, and the  ideals
$P_1,\dots,P_n$ are the maximal ideals in  $\mathbf W$.  Then there
exists an integral domain  $\mathbf E$  with an ideal  $J$  such that:

\noindent
{\bf{(\ref{theox}.1)}}
$\mathbf E$ $=$ $\mathbf W[x]'$
is a semi-local Dedekind domain that is
a finite integral extension domain of $\mathbf W$, where  $x$ $=$
$u^{\frac{1}{e}}$ in a fixed algebraic closure  $F(u)^*$  of  $F(u)$.

\noindent
{\bf{(\ref{theox}.2)}}
$[\mathbf E_{(0)}: \mathbf W_{(0)}]$ $=$ $e$.

\noindent
{\bf{(\ref{theox}.3)}}
$J$ $=$ $x \mathbf E$ is the Jacobson radical of  $\mathbf E$, 
$\mathbf E$  has
exactly  $n$  maximal ideals  $\mathbb M_1,\dots,\mathbb M_n$,
and  $\mathbb M_j \cap \mathbf W$ $=$ $P_j$  for  $j$
$=$ $1,\dots,n$.

\noindent
{\bf{(\ref{theox}.4)}}
$[(\mathbf E/\mathbb M_j):(\mathbf W/P_j)]$ $=$ $e_j$  for  $j$
$=$ $1,\dots,n$.

\noindent
{\bf{(\ref{theox}.5)}}
$u\mathbf E$ $=$ $J^e$.

\noindent
{\bf{(\ref{theox}.6)}}
The Rees integers of  $I\mathbf E$ $=$ $u\mathbf E$  are all equal to 
$e$.

\noindent
{\bf{(\ref{theox}.7)}}
If  $e_1 = \cdots = e_n = e$ and if there exists $b$ $\in$ $I$ 
such
that  $bV$ $=$ $IV$  for each Rees valuation ring  $(V,N)$  of  $I$,
then there exists a unit  $\theta$ $\in$  $\mathbf E$  such that
$\mathbf E$ $=$ $\mathbf W[\theta]$, so  $\mathbf E$  is a simple
free integral extension domain of  $\mathbf W$  of rank  $e$.

\noindent
{\bf{(\ref{theox}.8)}}
$\mathbf E$ $=$ $\mathbf T_{S'}$, where  $\mathbf T$ $=$
$R[u,tI,u^{\frac{1}{e}}]' \cap R[u^{\frac{1}{e}},t^{\frac{1}{e}}]$ and
${ S'}$ $=$
$\mathbf T - \cup \{q \mid q$  is a \textup{(}height one\textup{)}
associated prime ideal of  $u^{\frac{1}{e}} \mathbf T\}$.

\end{theo}

\begin{proof}
It is shown in Corollary \ref{co1y}(1) that:
$\mathbf W$  is a semi-local Dedekind
domain; $\mathbf Q$ is its Jacobson radical, and, the  ideals
$P_1,\dots,P_n$ are the maximal ideals in  $\mathbf W$.

(\ref{theox}.1) follows from Corollary \ref{co1y}(1)
and Corollary \ref{co1y}(3.2).

(Note:  $\mathbf E$  is frequently denoted  by $\mathbf D'$  in
Corollary \ref{co1y} and its proof.)

(\ref{theox}.2) follows from Corollary \ref{co1y}(2).

(\ref{theox}.3) is proved in Corollary \ref{co1y}(3.3) and Corollary
\ref{co1y}(3.2).

For (\ref{theox}.4), it is shown in Corollary \ref{co1y}(3) and
\ref{co1y}(3.2) that
there exist   $\theta_1,\dots,\theta_n$ $\in$ $\mathbf E$  such
that  $\mathbf E$ $=$ $\mathbf W[x,\theta_1,\dots,\theta_n]$ and
$\mathbb M_j$ $=$
$(P_j,x,{\theta_1},$ 
$\dots,{\theta_{j-1}},{\theta_{j+1}},\dots,{\theta_n}) \mathbf E$.

Therefore  $\mathbf E/\mathbb M_j$ $\cong$ $\mathbf W[\theta_j]/P_j$ 
$\cong$
$(W_j/Q_j)[ \theta_j + P_j]$,
and ${V_j}^*/{N_j}^*$  $=$  $W_j[x,\theta_j]/{N_j}^*$ $\cong$ 
$(W_j/Q_j)[\theta_j + P_j]$,
so  $\mathbf E/\mathbb M_j$ $\cong$ ${V_j}^*/{N_j}^*$ for  $j$ $=$ 
$1,\dots,n$.

Also,  $\mathbf W/P_j$ $\cong$ $W_j/Q_j$
for  $j$ $=$ $1,\dots,n$.
Therefore (\ref{theox}.4) follows from Theorem \ref{theo1y}.5.

Since  $I \mathbf E$ $=$ $u \mathbf E$,
(\ref{theox}.5) follows from Corollary \ref{co1y}(3.3), and then
(\ref{theox}.6) follows from (\ref{theox}.5) and Definition 
\ref{defi1}.5.
(\ref{theox}.7) is proved in Corollary \ref{co1y}(4).

Finally, for (\ref{theox}.8), $\mathbf E$ $=$ $\mathbf W[x]'$  is the
intersection of the Itoh (e)-valuation rings of  $I$,
by Corollary \ref{co1y}(3.1).  Also,  $\mathbf T_{S'}$
is the intersection of the Itoh (e)-valuation rings of  $I$,
by Definition \ref{DD}, so  $\mathbf E$ $=$ $\mathbf T_{S'}$, so
(\ref{theox}.8) holds.
\end{proof}

We next consider a powerful classical theorem of Krull, and to
state the theorem, we use the following terminology of Gilmer in
\cite{Gilmer}.

\begin{defi}
\label{consist}
{\em
Let $(V_1,N_1), \ldots, (V_n,N_n)$ be distinct
DVRs of a field $F$ and for  $j$ $=$
$1,\dots,n$,  let $K_j$ $=$ $V_j/N_j$ denote the residue field of
$V_j$. Let $m$ be a positive integer. By an {\bf m-consistent system
for} $\{ V_1, \ldots, V_n \}$, we mean a collection of sets $\mathbf S$ 
=
$\{ S_1, \ldots, S_n \}$ satisfying the following conditions:

\noindent
{\bf (1)}
For each  $j$, $S_j$ = $\{ (K_{j,i}, f_{j,i},e_{j,i}) \mid i = 1, 
\dots, s_j
\}$, where $K_{j,i}$ is a simple algebraic field extension of $K_j$
with $f_{j,i} = [K_{j,i}:K_j]$,  and  $s_j, e_{j,i} \in \mathbb N_{> 
0}$.

\noindent
{\bf (2)}
For each $j$, the sum $\sum_{i=1}^{s_j} e_{j,i}f_{j,i}$
= $m$.
}

\end{defi}

\begin{defi}
\label{real} {\em The  $m$-consistent system $\mathbf S$ as in
Definition~\ref{consist} is said to be {\bf realizable} if  there
exists a separable algebraic extension field $L$ of $F$ such that:

\noindent
{\bf (1)}
$[L : F]$ = $m$.

\noindent
{\bf (2)}
For $j$ $=$ $1,\dots,n$, $V_j$ has exactly $s_j$ extensions
$V_{j,1}, \ldots, V_{j,s_j}$ to $L$.

\noindent
{\bf (3)}
For $j$ $=$ $1,\dots,n$, the residue field of $V_{j,i}$ is
$K_j$-isomorphic to $K_{j,i}$,
and the ramification index of $V_{j,i}$ relative to $V_j$ is equal to  
$e_{j,i}$
(so $N_jV_{j,i}$ $=$ ${N_{j,i}}^{e_{j,i}}$).

\noindent If  $\mathbf S$  and  $L$  are as above, we say the field  $L$
{\bf{realizes}} $\mathbf S$  or that  $L$  is a {\bf{realization}} of
$\mathbf S$.
}

\end{defi}

\begin{rema}
\label{EXP}
{\bf{(\ref{EXP}.1)}}
{\em
With Definition \ref{real} in mind, and with the notation
of Theorem \ref{theox}, the semi-local Dedekind domain
$\mathbf E$ $=$ $\mathbf W[x]'$ (or its quotient field  $F(x)$)
in Theorem \ref{theox} is, except for separability,
a realization of the  $e$-consistent system $\mathbf S$ $=$
$\{S_1,\dots,S_n\}$  for the Rees valuation rings
$(W_1,Q_1),\dots,(W_n,Q_n)$  of  $uR[u,tI]$.
Here, $e$  is an arbitrary  positive common multiple
of the Rees integers  $e_1,\dots,e_n$  of  $uR[u,tI]$
(and of  $I$),
and for  $j$ $=$ $1,\dots,n$,
$S_j$ $=$ $\{(W_j/Q_j)[\overline{\theta_j}],e_j,\frac{e}{e_j})\}$,
where  $\theta_j$  is a root of
$X^{e_j} - w_j$  (with  $w_j$  playing the role of
$w$  in (\ref{co1y}.4.3) in the proof of Corollary \ref{co1y}),
and  $\overline \theta_j$ $=$ $\theta_j + Q_jW_j[\theta_j]$.

\noindent
{\bf{(\ref{EXP}.2)}}
More generally, let  $I$  be a nonzero proper ideal in a Noetherian
integral domain  $R$, let  $\mathbf R$ $=$ $R[u,tI]$  be the Rees
ring of  $R$  with respect to  $I$, let  $(W_1,Q_1), \ldots, (W_n,Q_n)$
be the  Rees valuation rings of $u\mathbf R$,  let  $\mathbf W$ $=$
$W_1 \cap \cdots \cap W_n$, and let  $k$ $\ge$ $2$  be an arbitrary
integer.  For  $j$ $=$ $1,\dots,n$, let  $e_j$  be the Rees integer of
$u \mathbf R$  with respect to  $W_j$, and let  $d_j$  be the greatest
common divisor of  $k$  and  $e_j$.
Also, let  $F$  be the quotient field of  $R$, let  $F(u)^*$  be
an algebraic closure of  $F(u)$, and let  $\mathbf E$  be the
integral closure of  $\mathbf W[x_k]$, where
$x_k$ $=$ $u^{\frac{1}{k}}$ $\in$ $F(u)^*$.  Then  $\mathbf E$
is a realization of the  $k$-consistent system $\mathbf S'$ $=$
$\{{S_1}',\dots,{S_n}'\}$  for the valuation rings
$(W_1,Q_1),\dots,(W_n,Q_n)$.
Here, for  $j$ $=$ $1,\dots,n$,
${S_j}'$ $=$ 
$\{(W_j/Q_j)[\overline{\theta_{j,k}}],d_j,\frac{k}{d_j})\}$,
where  $\theta_{j,k}$  is a root of
$X^{e_{j,k}} - w_{j,k}$  (with  $w_{j,k}$  playing the roll of
$w$  in (\ref{co1y}.4.3) in the proof of Corollary \ref{co1y}),
and  $\overline \theta_{j,k}$ $=$ $\theta_{j,k} + 
Q_jW_j[\theta_{j,k}]$.
}

\end{rema}

\begin{proof}
Item (\ref{EXP}.1) follows immediately from Theorem \ref{theox}
(it follows from Proposition \ref{prop6}(1) that
$Q_jW_j[\theta_j]$  is a prime ideal), and Item (\ref{EXP}.2) follows 
from
Proposition \ref{prop6}(3).
\end{proof}

\begin{theo}
\label{GK}
{\em (Krull \cite{Krull}):}
Let $(V_1,N_1), \ldots, (V_n,N_n)$
be distinct DVRs of a field $F$
with $K_j$ $=$ $V_j/N_j$
for  $j$ $=$
$1,\dots,n$, let  $m$ be a positive integer, and let
$\mathbf S$ = $\{ S_1, \ldots, S_n \}$ be an $m$-consistent system for
$\{ V_1, \ldots, V_n \}$
with
$S_j$ = $\{ (K_{j,i}, f_{j,i},e_{j,i}) \mid i = 1, \dots, s_j \}$
for  $j$ $=$ $1,\dots,n$.
Then $\mathbf S$ is realizable if one of the following
conditions is satisfied:

\begin{enumerate}[$(1)$]

\item
$s_j$ = $1$ for at least one $j$.

\item
$F$ has at least one DVR  $V$ distinct from
$V_1, \ldots, V_n $.

\item
For each monic polynomial $X^t + a_1X^{t-1} + \cdots + a_t$
with $a_i \in \cap_{j=1}^n V_j$ = $\mathbf D$,
and for each $h \in \mathbb N_{> 0}$
there exists an irreducible separable polynomial
$X^t + b_1X^{t-1} + \cdots + b_t \in \mathbf D[X]$ with
$b_l - a_l \in {N_j}^h$ for
each $l$ = $1, \ldots, t$ and  $j$ = $1, \ldots, n$.
\end{enumerate}

\end{theo}

\begin{observation*}~
{\em
\begin{enumerate}[$(a)$]

\item
Condition (1) of Theorem \ref{GK} is a property of the
$m$-consistent system $\mathbf S$ = $\{ S_1, \ldots, S_n \}$.

\item
Condition (2) of Theorem \ref{GK} is a property of the family of DVRs of
the field $F$.

\item
Condition (3) of Theorem \ref{GK} is a property of  the family 
$(V_1,N_1),
\ldots, (V_n,N_n)$.
\end{enumerate}
}
\end{observation*}

\begin{rema}
\label{rema3}
{\em
Let  $R$  be a Noetherian integral domain with quotient field  $F$,
and assume that $\operatorname{altitude}(R)$ $\ge$ $2$.
Then there exist infinitely many height one prime ideals in  $R$,
so  $R'$  has infinitely many height one prime ideals, by the
Lying-Over Theorem (\cite[(10.8)]{N2}).  Therefore, since  $R'$
is a Krull domain, by (\cite[(33.10)]{N2}),
there exist infinitely many distinct DVRs with quotient
field  $F$, hence
Theorem \ref{GK}(2) is always satisfied for such fields  $F$.
}

\end{rema}

We can now state and prove the first new result in this section.
It is closely related to Theorem \ref{theox}, and it is
also considerably more general.

\begin{theo}
\label{alternate}
Let  $R$  be a Noetherian integral domain with quotient field  $F$,
and assume that $\operatorname{altitude}(R)$ $\ge$ $2$.  Let  $I$
be a nonzero
proper ideal in  $R$, let
$(V_1,N_1),\dots,(V_{n},N_{n})$  \textup{(}$n$ $\ge$ $2$\textup{)} be 
the  Rees
valuation rings of  $I$, let  $\mathbf D$ $=$ $V_1 \cap \cdots \cap 
V_n$
 \textup{(}so
$\mathbf D$  is a Dedekind domain with exactly  $n$  maximal ideals  
$M_j$
$=$
$N_j \cap \mathbf D$  \textup{(}$j$ $=$ $1,\dots,n$\textup{))}, and 
let
${M_1}^{e_1} \cdots {M_n}^{e_n}$
\textup{(}$=$ ${M_1}^{e_1} \cap \cdots \cap {M_n}^{e_n}$\textup{)}
be an irredundant primary
decomposition of  $I\mathbf D$  \textup{(}so  $e_1,\dots,e_n$  are the 
Rees
integers
of  $I$\textup{)}.  Let  $m$ be the least common multiple of
$e_1,\dots,e_n$, let  $d_j$ $=$ $\frac{m}{e_j}$ \textup{(}$j$ $=$
$1,\dots,n$\textup{)},
and  let  $e$  $=$ $km$  be a positive multiple of  $m$.  Then there
exists an integral domain  $\mathbf E$  with an ideal  $J$  such that:

\noindent
{\bf{(\ref{alternate}.1)}}
$\mathbf E$ $=$ $\mathbf D[\theta]$ is a semi-local Dedekind domain 
that
is a simple free separable integral extension
domain of  $\mathbf D$.

\noindent
{\bf{(\ref{alternate}.2)}}
$[\mathbf E:\mathbf D]$ $=$ $e$.

\noindent
{\bf{(\ref{alternate}.3)}}
$J$  is the Jacobson radical of  $\mathbf E$, and  $\mathbf E$  has 
exactly
  $k(e_1 + \cdots + e_n)$  maximal ideals.

\noindent
{\bf{(\ref{alternate}.4)}}
$[(\mathbf E/Q_{j,i}) : (\mathbf D/M_j)]$ $=$ $1$  for
$j$ $=$ $1,\dots, n$  and  for the  $ke_j$   associated
prime ideals  $Q_{j,i}$  of  $M_j\mathbf E$.

\noindent
{\bf{(\ref{alternate}.5)}}
$I\mathbf E$ $=$ $J^m$.

\noindent
{\bf{(\ref{alternate}.6)}}
The Rees integers of  $I\mathbf E$  are all equal to  $m$.

\end{theo}

\begin{proof}
Let   $\mathbf S$ $=$ $\{S_1,\dots,S_n\}$ with $S_j$
= $\{(K_{j,i},1,d_j) \mid i = 1,\dots,ke_j\}$
($j$ $=$ $1,\dots,n$).
Observe that $d_j \cdot ke_j$ $=$ $km$ $=$ $e$  for  $j$
$=$ $1,\dots,n$, so  $\mathbf S$  is an $e$-consistent system
for  $\{\mathbf D_{M_1},\dots,\mathbf D_{M_n}\}$, so  $\mathbf S$
is a realizable  $e$-consistent system for
$\{\mathbf D_{M_1},\dots,\mathbf D_{M_n}\}$,
by Theorem \ref{GK}(2) and Remark \ref{rema3}.
Therefore the integral closure
$\mathbf E$  of  $\mathbf D$  in a realization  $L$  of  $\mathbf S$ for
$\{\mathbf D_{M_1},\dots,\mathbf D_{M_n}\}$  is a simple free separable
integral extension domain
of  $\mathbf D$  such that  $[\mathbf E:\mathbf D]$ $=$ $e$  and  
$\mathbf E$
is a Dedekind domain (by
\cite[Theorem 19, p. 281]{ZS1}).  Also, since $S_j$ = 
$\{(K_{j,i},1,d_j) \mid i = 1,\dots,ke_j \}$
($j$ $=$ $1,\dots,n$):
each  $V_j$  has exactly  $ke_j$  extensions  $(V_{j,i},N_{j,i})$  to  
$L$;
$\mathbf E$ $=$ $V_{1,1} \cap \cdots \cap V_{n,ke_n}$  and  $Q_{j,i}$ 
$
=$
$N_{j,i} \cap \mathbf E$; $\mathbf E/Q_{j,i}$ $\cong$
$V_{j,i}/N_{j,i}$ $\cong$ $V_j/N_j$ $\cong$ $\mathbf D/M_j$  ($j$ $=$ 
$1,
\dots,n$  and  $i$
$=$ $1,\dots,ke_j$);
and, $N_jV_{j,i}$ $=$ ${N_{j,i}}^{d_j}$ ($j$ $=$ $1,\dots,n$  and
$i$ $=$ $1,\dots,ke_j$).

Therefore (\ref{alternate}.1) - (\ref{alternate}.4) hold, and
$$I\mathbf E = \prod_{j=1}^n ({M_j}^{e_j}\mathbf E) = 
\prod_{j=1}^n
  (\prod_{i=1}^{ke_j} ({N_{j,i}}^{d_j}))^ {e_j} = (\prod_{j=1}^n
\prod_{i=1}^{ke_j}  {N_{j,i}})^m = J^m,$$
where  ${J}$ = $\prod_{j=1}^n  \prod_{i=1}^{ke_j}  {N_{j,i}}$
is the Jacobson radical  $J$  of  $\mathbf E$.
Thus  $I\mathbf E$ $=$ $J^e$, so (\ref{alternate}.5) holds, and since
$\mathbf E$  is a semi-local domain with
Jacobson radical  $J$, (\ref{alternate}.6) follows immediately from
(\ref{alternate}.5),
Definition \ref{defi1}.5, and Remark \ref{RR}.4.
\end{proof}

\begin{rema}
\label{extra0}
{\em
In the proof of Proposition \ref{alternate},
there are many cases when the simpler  $e$-consistent system
$\mathbf T$ $=$ $\{T_1,\dots,T_n\}$ with $T_j$
= $\{(K_{j,1},ke_j,d_j)\}$  ($j$ $=$ $1,\dots,n$)  can be used
in place of  $S_j$
= $\{(K_{j,i},1,d_j) \mid i = 1,\dots,ke_j\}$
($j$ $=$ $1,\dots,n$),
and then the resulting realization  $\mathbf E$  of  $\mathbf T$
for  $\{\mathbf D_{M_1},\dots,\mathbf D_{M_n}\}$  has the same number of
maximal ideals as  $\mathbf D$.  However, $\mathbf T$  cannot
always be used, since, for example, if  $\mathbf D/M_j$  is
algebraically closed and  $ke_j$ $\ge$ $2$, then there are no extension
fields  $K_{j,1}$
of  $K_j$ $=$ $\mathbf D/M_j$  such that  $[K_{j,1}:K_j]$ $=$ 
$ke_j$.
}

\end{rema}

\begin{rema}
\label{extra}
{\em
In the proof of Proposition \ref{alternate}, if   $S_j$
= $\{(K_{j,i},1,d_j) \mid i = 1,\dots,ke_j\}$  is replaced with  
$U_j$
$=$ $\{(K_{j,i},1,kd_j) \mid i = 1,\dots,e_j \}$ ($j$ $=$ 
$1,\dots,n$),
then the same conclusions hold, but replace:

\noindent
{\bf (a)}
``each  $V_j$  has exactly  $ke_j$  extensions  $(V_{j,i},N_{j,i})$  to
$L $'' with
``each  $V_j$  has exactly  $e_j$  extensions  $(V_{j,i},N_{j,i})$ to 
$L$'';

\noindent
{\bf (b)}
``the Rees integers of  $I\mathbf E$  are all equal to  $m$,'' with
``the Rees integers of  $I\mathbf E$  are all equal to  $e$'';

\noindent
{\bf (c)}
``$I\mathbf E$ $=$ $J^m$''  with  ``$I\mathbf E$ $=$ $J^e$'';

\noindent
{\bf (d)}
``$N_jV_{j,i}$ $=$ ${N_{j,i}}^{d_j}$ ($j$ $=$ $1,\dots,n$  and $i$
$= $ $1,\dots,ke_j$)''
with
``$N_jV_{j,i}$ $=$ ${N_{j,i}}^{kd_j}$ ($j$ $=$ $1,\dots,n$  and $i$ =
$1,\dots,e_j$)'';

\noindent
{\bf (e)}
``$I\mathbf E$ = $\prod_{j=1}^n ({M_j}^{e_j}\mathbf E)$ =
$\prod_{j=1}^n  (\prod_{i=1}^{ke_j} ({N_{j,i}}^{d_j}))^{e_j}$ =
$(\prod_{j=1}^n  \prod_{i=1}^{ke_j}  {N_{j,i}})^m$''
with
``$I \mathbf E $  = \, $\prod_{j=1}^n ({M_j}^{e_j}\mathbf E)$ =
$\prod_{j=1}^n  (\prod_{i=1}^{e_j} ({N_{j,i}}^{kd_j}))^ {e_j} $ 
$=$
$(\prod_{j=1}^n  \prod_{i=1}^{e_j}  {N_{j,i}})^e$''.
}

\end{rema}

Our first corollary of Theorem \ref{alternate} is a more
complete and detailed version of \cite[Theorem 2.8(1)]{HRR1x}.

\begin{coro}
\label{co1}
Let  $I$  be a nonzero proper ideal in a Noetherian domain  $R$, let
$(V_1,N_1),$ $\dots,$ $(V_{n},N_{n})$  \textup{(}$n$ $\ge$ $2$\textup{)} 
be
the  Rees
valuation rings of  $I$, let
$e_j$  be the Rees integer of  $I$  with respect to  $V_j$  
\textup{(}$j$
$=$ $1,\dots,n$\textup{)}, and let  $e$ $=$ $km$
be a positive multiple of the least common multiple $m$ of
$e_1, \dots, e_n$.  Then there exists an integral domain  $B_e$
such that:

\noindent
{\bf (1)}
$B_e$  is a semi-local Dedekind domain that is a
simple free separable integral extension domain of  $R$.

\noindent
{\bf (2)}
$[B_e:R]$ $=$ $e$.

\noindent
{\bf (3)}
The Rees integers of  $IB_e$  are all equal to  $m$.
\end{coro}

\begin{proof}
Let  $\mathbf D$  be the intersection of the  $n$  Rees valuation rings
$(V_j,N_j)$
of  $I$, and let  $M_j$ $=$ $N_j \cap \mathbf D$  ($j$ $=$ 
$1,\dots,n$), so
$\mathbf D$  is a semi-local Dedekind domain with exactly  $n$  maximal
ideals  $M_j$
and  $I\mathbf D$ $=$ ${M_1}^{e_1} \cap \cdots \cap {M_n}^{e_n}$
$=$ ${M_1}^{e_1} \cdots {M_n}^{e_n}$.
Therefore by Proposition \ref{alternate} there exists
a simple free separable
extension field  $L_e$  of the quotient field  $F$  of  $R$
such that:  $[L_e:F]$ $=$ $e$;
the integral closure  $\mathbf E_e$  of $\mathbf D$ in $L_e$ is a
finite (by
free separability) integral
extension domain of  $\mathbf D$  and is a semi-local Dedekind
domain with exactly  $ke_j$  maximal ideals  $Q_{j,i}$  lying over
each  $M_j$  ($j$ $=$ $1,\dots,n$); and,
$(I\mathbf D)\mathbf E_e$ $=$ ${J_e}^m$, where  $J_e$
is the Jacobson radical of  $\mathbf E_e$.
By separability  $L_e$ $=$ $F[\theta]$, so there exists  $r$
$\in$ $R$  such that
$r \cdot \theta$  is integral over  $R$, so let  $B_e$ $=$
$R[r \cdot \theta]$.  Then  $B_e$  is a simple free separable integral
extension domain of  $R$, $[B_e:R]$ $=$ $e$, and  $\mathbf E_e$  is 
the
intersection of the Rees valuation rings of  $IB_e$
by Remark \ref{RR}.6.  Since  $(IB_e)\mathbf E_e$ $=$
$(I\mathbf D)\mathbf E_e$ $=$ ${J_e}^m$ and
since  $J_e$  is the Jacobson radical of  $\mathbf E_e$, it follows
that  the Rees integers of  $IB_e$  are all equal to  $m$.
\end{proof}

\begin{rema}
\label{extra2}
{\em

\noindent
{\bf{(\ref{extra2}.1)}}
With the notation of Corollary \ref{co1},
there exists an integral domain  $C_e$  such that:

\noindent
{\bf (1)}
$C_e$  is a semi-local Dedekind domain that is a
simple free separable integral extension domain  of  $R$.

\noindent
{\bf (2)}
$[C_e:R]$ $=$ $e$.

\noindent
{\bf (3)}
The Rees integers of  $IC_e$  are all equal to  $e$.

\noindent
{\bf{(\ref{extra2}.2)}}
Let  $I$  be a regular proper ideal in a Noetherian ring  $R$, let
$(V_1,N_1),\dots,(V_{n},N_{n})$  ($n$ $\ge$ $2$) be the  Rees
valuation rings of  $I$, let
$e_j$  be the Rees integer of  $I$  with respect to  $V_j$  ($j$
$=$ $1,\dots,n$), and let  $e$ $=$ $km$
be a positive multiple of the least common multiple $m$ of
$e_1, \dots, e_n$.  Assume that $\operatorname{Rad}(0_R)$ is prime,
say $\operatorname{Rad}(0_R)$ $=$ $z$.
Then there exists a ring  $B_e$  such that:

\noindent
{\bf (1)}
$B_e$  is a simple free integral extension ring  of  $R$.

\noindent
{\bf (2)}
$[B_e:R]$ $=$ $e$.

\noindent
{\bf (3)}
The Rees integers of  $IB_e$  are all equal to  $m$.

\noindent
{\bf (4)}
$zB_e$ $=$ $\operatorname{Rad}(0_{B_e})$, so  $zB_e$
is the only minimal prime ideal in  $B_e$.

\noindent
{\bf{(\ref{extra2}.3)}}
With the notation of (\ref{extra2}.2),
there exists a ring  $C_e$  such that:

\noindent
{\bf (1)}
$C_e$  is a simple free integral extension ring  of  $R$.

\noindent
{\bf (2)}
$[C_e:R]$ $=$ $e$.

\noindent
{\bf (3)}
The Rees integers of  $IC_e$  are all equal to  $e$.

\noindent
{\bf (4)}
$zC_e$ $=$ $\operatorname{Rad}(0_{C_e})$, so  $zC_e$
is the only minimal prime ideal in  $C_e$.
}

\end{rema}

\begin{proof}
The proof of (\ref{extra2}.1) is the same as the proof of
Corollary \ref{co1}, but use
Remark \ref{extra} in place of Proposition \ref{alternate}.

For (\ref{extra2}.2), since  $\operatorname{Rad}(0_R)$ $=$ $z$ is
the only minimal prime ideal in  $R$,
the Rees valuation rings of  $I$  are the Rees valuation rings of
$\overline I$ $=$ $(I+z)/z$ (see Remark \ref{RR}.3).  Also, by 
Corollary
\ref{co1} there exists an
integral domain  $\overline{B_e}$  such that:

\noindent
(1')
$\overline{B_e}$  is a semi-local Dedekind domain that is a
simple free separable integral extension domain  of  $\overline R$ $=$
$R/z$.

\noindent
(2')
$[\overline{B_e}:\overline R]$ $=$ $e$.

\noindent
(3')
The Rees integers of  $\overline I \overline{B_e}$  are all
equal to  $m$.

\noindent
Let  $f_e(X)$  be the pre-image in  $R[X]$  of the irreducible monic
polynomial  $\overline {f_e}(X)$  (of degree  $e$)
in $\overline R[X]$  such that  $\overline {B_e}$ $=$
$\overline R[X]/(\overline {f_e}(X) \overline R[X])$.
Then  $\overline {B_e}$ $=$ $R[X]/((f_e(X),z)R[X])$ $=$ 
$R[x]/(zR[x])$,
where  $x$ $=$ $X+(f_e(X)R[X])$.

Let  ${B_e}$ $=$ $R[x]$.  Then, since  $\overline {B_e}$ $=$
$R[x]/(zR[ x])$  is an integral
domain, it follows that  $zB_e$
is a prime ideal.  Also, by hypothesis, there exists
$r$ $\in$ $\mathbb N_{> 0}$
such that  $z^r$ $=$ $(0)$  in  $R$, so it follows that  $zB_e$  is 
the
only minimal
prime ideal in  $B_e$, so (4) holds.  Therefore the Rees
valuation rings and Rees integers of  $IB_e$  are the Rees
valuation rings and Rees integers of  $\overline I \overline{B_e}$,
by Remark \ref{RR}.3, so (1) - (3) follow immediately from (1') - (3').

The proof of (\ref{extra2}.3) is the same as the proof of 
(\ref{extra2}.2),
but use
Remark \ref{extra2}.1 in place of Corollary \ref{co1}.
\end{proof}

The next corollary of Theorem \ref{alternate} extends
Corollary \ref{co1} to Noetherian rings.

\begin{coro}
\label{co2}
Let  $I$  be a regular proper ideal in a Noetherian ring $R$, let
$e_1,\dots,e_n$  be the Rees integers of  $I$, let  $m$
be the least common multiple  of $e_1,\dots, e_n$, and let
$e$ $=$ $km$  be a positive multiple of  $m$.
Then there exist rings  $R^*$  and  $B_e$  such that:

\noindent
{\bf (1)}
$R$ $\subseteq$ $R^*$ $\subseteq$ $B_e$ and $B_e$ is a
finite integral
extension ring  of  $R$.

\noindent
{\bf (2)}
$IR^*$  and  $IB_e$  are  regular proper ideals.

\noindent
{\bf (3)}
There is a one-to-one correspondence between the minimal
prime ideals  $w^*$
in  $B_e$  such that  $IB_e +w^*$ $\ne$ $B_e$,
the  minimal prime ideals  $z^*$ in $R^*$ such that $IR^*+z^*$ $\ne$
$R^*$, and
the  minimal prime ideals $z$ in $R$ such that $I+z$ $\ne$ $R$; namely,
$w^*$ $=$ $z^*B_e$  and  $z$ $=$ $z^* \cap R$.

\noindent
{\bf (4)}
The Rees integers of  $IB_e$  are all equal to  $e$.

\end{coro}

\begin{proof}
Let  $(0)$ $=$ $\cap \{q_h \mid h = 1,\dots,k\}$  be an irredundant
primary decomposition of the zero ideal in  $R$  and
let  $\operatorname{Rad}(q_h)$ $=$ $z_h$
($h$ $=$ $1,\dots,k$).  Assume that the  $z_h$  are re-ordered
so that  $z_1,\dots,z_{d_1}$  are the minimal
prime ideals $z$  in  $R$  such that
$I+z$ $\ne$ $R$  and
$z_{d_{1+1}},\dots,z_{d_2}$  are the remaining minimal
prime ideals in  $R$  (so  $d_1$ $\le$ $d_2$ $\le$ $k$).

Rewrite $\cap \{q_h \mid h = 1,\dots,k\}$  as
$\cap \{Z_h \mid h = 1,\dots,d_2\}$, where  $Z_h$  is the
intersection of all  $q_i$  such that  $z_i$ $\supseteq z_h$,
for  $h$ $=$ $1,\dots,d_2$.
For $i$ $=$ $1,\dots,d_2$, let $R_i$ $=$ $R/Z_i$,
let  $\overline{z_i}$ $=$ $z_i/Z_i$  
(so the unique minimal prime ideal in  $R_i$ is  $\overline{z_i}$,
therefore,  $\overline{z_i}$ $=$ $\operatorname{Rad}(0_{R_i}))$, and
let  $I_i$ $=$ $(I+Z_i)/Z_i$  (so  $I_i$  is a regular proper
ideal in  $R_i$ for  $i$ $=$ $1,\dots,d_1$  and  $I_i$ $=$
$R_i$  for  $i$ $=$ $d_1 +1,\dots,d_2$).

By Remark \ref{RR}.3 the set of Rees valuation rings
of  $I$  is the disjoint union of the sets of Rees valuation
rings of the ideals $I_i$  ($i$ $=$ $1,\dots,d_1$), so the
Rees integers of the  $I_i$  are among the Rees
integers of  $I$. 
Hence $e$  is a multiple of the least common multiple of
the Rees integers of the ideals  $I_i$  ($i$ $=$ $1,\dots,d_1$).

Therefore, by Remark \ref{extra2}.3,
for $i$ $=$ $1,\dots,d_1$,  there
exists a simple free integral extension ring  $B_{i,e}$  of  $R_i$
such that
$[B_{i,e}:R_i]$ $=$ $e$,
the Rees integers of  $I_i B_{i,e}$  are all equal to  $e$,
and $\overline{z_i}B_{i,e}$  is the only minimal prime ideal in  
$B_{i,e}$.

Let  $R^*$ $=$ $R_1 \bigoplus \cdots \bigoplus R_{d_2}$,
let  $I^*$ $=$ $I_1$$ \bigoplus \cdots \bigoplus$$ I_{d_2}$
so  $I^*$ $=$ $IR^*$, and for $i$ $=$ $1,\dots,d_2$,  let
${z_i}^*$ $=$ $R_1$ $\bigoplus \cdots \bigoplus$ $R_{i-1}$
$\bigoplus$ $\overline{z_i}$ $\bigoplus$$ R_{i+1}$
$ \bigoplus \cdots \bigoplus$$ R_{d_2}$.
Then $R^*$  is a finite integral extension ring of  $R$,
$I^*$  is a regular proper ideal in $R^*$, and for  $i$ $=$
$1,\dots,d_1$,  the  ${z_i}^*$  are the 
minimal prime ideals  $z^*$  in  $R^*$  such that  $I^* + z^*$ $\ne$ $R^*$.
Also,  ${z_i}^* \cap R$ $=$ $z_i$  and
$R^*/{z_i}^*$ $=$ $R_i/\overline{z_i}$   $=$ $R/z_i$ ($i$
$=$ $1,\dots,d_2$),
so the Rees valuation rings of  $I^*$  are the
Rees valuation rings of  $I$,
so the ideals  $I$  and  $I^*$  have the same
Rees integers.

Let $B_e$ $=$ $B_{1,e}$$ \bigoplus \cdots \bigoplus$$ B_{d_2,e}$,
where, for notational convenience, we let  $B_{h,e}$ $=$ $R_h$  for  
$h$ $=$ $d_1 +1,\dots,d_2$.
Also, for  $i$ $=$ $1,\dots,d_2$,  let ${w_i}^*$ $=$
$B_{1,e}$$ \bigoplus \cdots \bigoplus $$B_{i-1,e}  $$\bigoplus 
$$\overline{z_i}B_{i,e} $$\bigoplus $$B_{i+1,e} $$\bigoplus \cdots \bigoplus 
$$B_{d_2,e}$,
so  ${w_i}^* \cap R^*$ $=$ ${z_i}^*$  and  ${w_i}^*$ $=$ 
${z_i}^*B_e$.
Then $B_e$  is a finite  $R^*$-module (and is also a
finite integral extension ring of  $R$), $IB_e$ $=$ $I^*B_e$ $=$
$IB_{1,e}$ $ \bigoplus \cdots \bigoplus$ $ IB_{d_1,e}$ $ \bigoplus$
$B_{d_1+1,e }$ $ \bigoplus \cdots \bigoplus$ $B_{d_2,e}$
is a regular proper ideal in  $B_e$,
the  ${w_i}^*$  are  $d_2$  minimal prime ideals  in  $B_e$, and
$B_e/{w_i}^*$ $=$ $B_{i,e}/(\overline{z_i}B_{i,e})$  for  $i$ $=$
$1,\dots,d_2$.

Since the  ideals  ${z_i}^*$
are the minimal prime ideals  $z^*$  in  $R^*$ such that $I^* + z^*$ 
$\ne
$ $R^*$  (for  $i$ $=$ $1,\dots,d_1$),
and since  ${w_i}^*$  is a minimal prime ideal in  $B_e$  such that
${w_i}^*$ $=$ ${z_i}^*B_e$,
it follows that the  ${w_i}^*$  ($i$ $=$ $1,\dots,d_1$)  are the
minimal prime ideals  $w^*$  in  $B_e$  such that
$IB_e + w^*$ $=$ $I^*B_e+w^*$ $\ne$ $B_e$. So the set of Rees valuation rings of  $I^*B_e$  is the disjoint union 
of the sets of Rees valuation rings of the ideals $(I^*B_e+{w_i}^*)/{w_i}^*$ 
$=$ $(I^*B_e+{z_i}^*B_e)/{w_i}^*$ $=$
$(B_{1,e}$$ \bigoplus \cdots \bigoplus$$ B_{i-1,e}$$ \bigoplus$$
(I_i+\overline{z_i})B_{i,e}$$ \bigoplus$$ B_{i+1,e}$$
\bigoplus \cdots \bigoplus$$ B_{d_2,e})/{w_i}^*$
$=$ $(I_iB_{i,e} + \overline{z_i}B_{i,e})/(\overline{z_i}B_{i,e})$.
Therefore,
since  $\overline{z_i}B_{i,e}$  is the unique minimal prime ideal in
$B_{i, e}$, it follows
that, for  $i$ $=$ $1,\dots,d_1$,  $I_iB_{i,e}$ and  $(I_iB_{i,e}$ $+$
$ \overline{z_i}B_{i,e})/(\overline{z_i}B_{i,e})$  have
the same Rees valuation rings and the same Rees integers.

Finally, for  $i$ $=$$1,\dots,d_1$,
the Rees integers of  $I_iB_{i,e}$  are all equal to  $e$
(by the last sentence in the second preceding paragraph),
so it follows that
the Rees integers of  $IB_e$  are all equal to  $e$.
\end{proof}

To state an additional corollary, we need the following definition.

\begin{defi}
\label{defi3}
{\em
Let  $I$  and  $J$  be ideals in a ring  $R$.  Then:

\noindent
{\bf{(\ref{defi3}.1)}}
$I$  and  $J$  are
{\bf{projectively equivalent}}
in case there exist  $i,j$ $\in$ $\mathbb N_{> 0}$  such
that  $(I^i)_a$ $=$ $(J^j)_a$
(see (\ref{defi1}.2)).

\noindent
{\bf{(\ref{defi3}.2)}}
$I$  is
{\bf{projectively full}}
in case, for each ideal  $J$  in  $R$  that is projectively
equivalent to  $I$  (see (\ref{defi3}.1)), $J_a$ $=$ $(I^k)_a$
for some  $k$ $\in$ $\mathbb N_{> 0}$.
}

\end{defi}

\begin{rema}
\label{xxxx}

{\em
With Definition \ref{defi3} in mind, it should be noted
that Theorem \ref{alternate}.6 shows that
the Jacobson radical  $J$  of  $\mathbf E$ is projectively equivalent to 
$I
\mathbf E$, and
since  $\mathbf E$  is a semi-local Dedekind domain, it follows
that  $J$  is a projectively full
radical ideal  whose Rees integers are all equal to one.
}

\end{rema}

The next corollary of Theorem \ref{alternate} was proved
in \cite[Theorem 2.8(2)]{HRR1x} in the case when  $R$  is a Noetherian
domain of altitude one.

\begin{coro}
\label{extend}
If  $R$  is a Noetherian ring of
altitude one, then for each regular
proper ideal  $I$  in  $R$  there exists
a finite integral extension ring  $A$  of
$R$  with an ideal  $J$  such that
$J$  and  $IA$  are projectively equivalent,
$J$  is a projectively full radical
ideal, and the Rees integers of  $I$  are
all equal to one.
\end{coro}

\begin{proof}
It is shown in \cite[Theorem 2.8(2)]{HRR1x} that
this result holds for nonzero proper ideals in
Noetherian domains of altitude one.  So a proof
similar to the proof of Corollary \ref{co2} shows
that it continues to hold for regular
proper ideals in Noetherian rings of altitude one.
\end{proof}

\vspace{.15in}
\begin{flushleft}
Department of Mathematics, University of California, Riverside,
California 92521-0135
{\em E-mail address: youngsu.kim@ucr.edu}

\vspace{.15in}

Department of Mathematics, University of California, Riverside,
California 92521-0135
{\em E-mail address: ratliff@math.ucr.edu}

\vspace{.15in}

Department of Mathematics, University of California, Riverside,
California 92521-0135
{\em E-mail address: rush@math.ucr.edu}

\end{flushleft}


\bigskip
\begin{thebibliography}{99}

\bibitem{CHRR3}
C. Ciuperc\v{a}, W. J. Heinzer, L. J. Ratliff, Jr., and D. E. Rush,
{\em Projectively full ideals in Noetherian rings \textup{(II)}}, J. Algebra
305 (2006), 974-992.

\bibitem{End}
O. Endler,
{\em Valuation Theory},
Springer-Verlag, New York, 1972.

\bibitem{Gilmer}
R. Gilmer, {\em Pr\"{u}fer domains and rings of integer-valued
polynomials}, J. Algebra, 129 (1990), 502-517.

\bibitem{HRR1x}
W. J. Heinzer, L. J. Ratliff, Jr., and D. E. Rush, {\em
Projective equivalence of ideals in Noetherian integral domains},
J. Algebra 320 (2008), 2349-2362.

\bibitem{It}
S.  Itoh,
{\em
Integral closures of ideals generated by regular sequences},
J. Algebra 117 (1988), 390-401.

\bibitem{Krull}
W. Krull, {\em \"{U}ber einen Existensatz der Bewertungstheorie},
Abh. Math. Sem. Univ. Hamburg, 23 (1959), 29-55.

\bibitem{N2}
M. Nagata,
{\em Local Rings},
Interscience, John Wiley, New York, 1962.

\bibitem{Rt}
L. J. Ratliff, Jr.,
{\em A Characterization of analytically unramified semi-local rings
and applications},
Pacific J. Math. 27 (1968), 127 - 143.

\bibitem{R2}
L. J. Ratliff, Jr.,
{\em On prime divisors of the integral closure of a principal ideal},
J. Reine Angew. Math. 255 (1972), 210 - 220.

\bibitem{ReF}
D. Rees,
{\em A note on form rings and ideals},
Mathematika 4 (1957), 51-60.

\bibitem{SH}
I. Swanson and C. Huneke,
{\em Integral Closure of Ideals, Rings and Modules},
Cambridge Univ. Press, Cambridge, 2006.

\bibitem{ZS1}
O. Zariski and P. Samuel,
{\em Commutative Algebra, Vol. 1},
D. Van Nostrand, New York, 1958.

\bibitem{ZS2}
O. Zariski and P. Samuel,
{\em Commutative Algebra, Vol. 2},
D. Van Nostrand, New York, 1960.

\end{thebibliography}
\end{document}